\documentclass[12pt,a4paper,final]{iopart} 

\usepackage{a4wide,amssymb}
\usepackage{todonotes}
\usepackage{iopams}
\usepackage{hyperref}
\usepackage{amsthm}
\usepackage{subcaption}
\usepackage{epstopdf}
\usepackage{multirow}
\graphicspath{{Images/}}
\usepackage{pgfplotstable, booktabs, longtable} 
\usepackage{float} 

\newcommand{\R}{\mathbb{R}}
\newtheorem{theorem}{Theorem}[section]

\newtheorem{lemma}[theorem]{Lemma}
\theoremstyle{definition}

\begin{document}

\title{Homoclinic saddle to saddle-focus transitions in 4D systems}
\author{Manu Kalia$^1$, Yuri A. Kuznetsov$^{1,2}$ and Hil G.E. Meijer$^1$}

\begin{abstract}
A saddle to saddle-focus homoclinic transition when the stable leading eigenspace is 3-dimensional (called the 3DL bifurcation) is analyzed. Here a pair of complex eigenvalues and a real eigenvalue exchange their position relative to the imaginary axis, giving rise to a 3-dimensional stable leading eigenspace. This transition is different from the standard Belyakov bifurcation, where a double real eigenvalue splits either into a pair of complex-conjugate eigenvalues or two distinct real eigenvalues. In the wild case, we obtain sets of codimension 1 and 2 bifurcation curves and points that asymptotically approach the 3DL bifurcation point and have a structure that differs from that the standard Belyakov case. We also give an example of this bifurcation in the wild case occuring in a perturbed Lorenz-Stenflo 4D ODE model.
\end{abstract}
\address{$^1$ Department of Applied Mathematics, University of Twente, Zilverling Building, P.O. Box 217, 7500AE Enschede, The Netherlands}
\address{$^2$ Mathematical Institute, Utrecht University, Budapestlaan 6, 3584CD Utrecht, The Netherlands}
\ead{[M.Kalia, I.A.Kouznetsov, H.G.E.Meijer]@utwente.nl}

\section{Introduction}

Homoclinic orbits play an important role in the analysis of ODEs depending on parameters
\begin{equation}
\dot{x}=F(x,\alpha), \quad x\in\R^n, \alpha\in\R^m, 
\label{eq:intro1}
\end{equation}
where $F$ is smooth. Orbits homoclinic to \textit{hyperbolic equilibria} are of specific interest, as they are structurally unstable, and the corresponding parameter values generically belong to codim 1 manifolds in the parameter space $\R^m$. Bifurcations in generic one-parameter families transverse to such manifolds depend crucially on the configuration of \textit{leading eigenvalues} of the equilibrium, i.e. the stable eigenvalues with largest real part, and the unstable eigenvalues with smallest real part.\\

\begin{figure}[htbp]
	\centering
	\begin{subfigure}[b]{0.25\textwidth}
		\includegraphics[width=\textwidth]{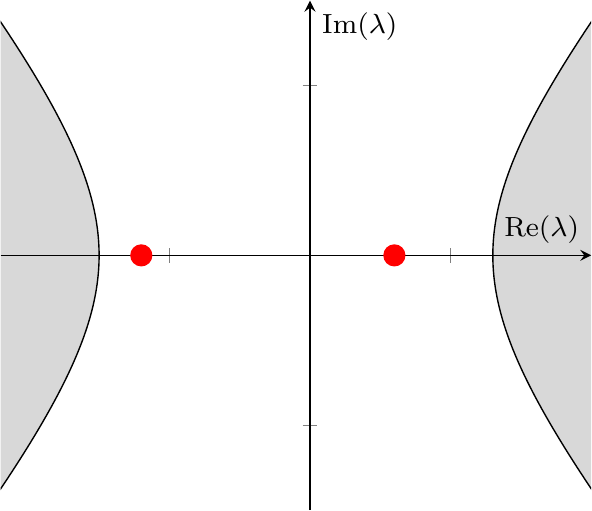}
		\caption{Saddle}
	\end{subfigure}\hfill
	\begin{subfigure}[b]{0.25\textwidth}
		\includegraphics[width=\textwidth]{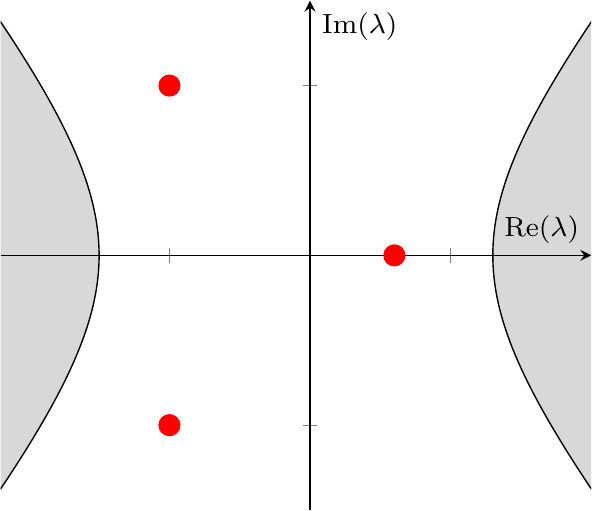}
		\caption{Saddle-focus}
	\end{subfigure}\hfill
	\begin{subfigure}[b]{0.25\textwidth}
		\includegraphics[width=\textwidth]{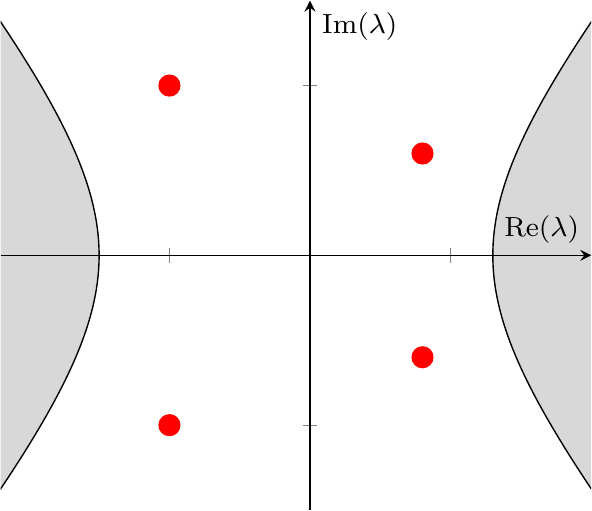}
		\caption{Focus-focus}
	\end{subfigure}
	\caption{Configurations of leading eigenvalues $\lambda$ (red). Gray area contains all non-leading eigenvalues. }
\label{fig:intro1}
\end{figure}

In \autoref{fig:intro1}, we see three configurations with simple leading eigenvalues, for which a detailed description of the bifurcations occurring near the homoclinic orbit is available (see, e.g. \cite{DS-5,IlyaLi1999,Shilnikovs2001,HomburgSandstede2010}). For example, in the saddle case, a single periodic orbit appears generically. In the saddle-focus case, we could assume that the leading stable eigenvalues are complex by applying time-reversal if necessary. In this case, infinitely many periodic orbits exist if the \textit{saddle quantity} $\sigma_0$, defined as the sum of the real part of the leading unstable and stable eigenvalues, is positive. This phenomenon is called \textit{Shilnikov's homoclinic chaos}
\cite{Shilnikov1965,Shilnikov2007}. On the contrary, if $\sigma_0$ is negative, then generically only one periodic orbit appears. Thus, the sign of $\sigma_0$ distinguishes \textit{wild} and \textit{tame} saddle-focus homoclinic cases. Note that in the wild case many other bifurcations occur nearby, including infinite sequences of fold (limit point, LP) and period-doubling (PD) bifurcations of periodic orbits, as well as secondary homoclinic bifurcations, which all accumulate on the primary homoclinic bifurcation manifold \cite{Gonchenko1997}. In the focus-focus case, which will not be considered in this paper, homoclinic chaos is always present.\\

Moving along the primary homoclinic manifold in the parameter space of (\ref{eq:intro1}), one may encounter a transition from the saddle case (a) to the saddle-focus case (b). This is a degenerate situation, and the corresponding homoclinic parameter values form generically a codim 2 sub-manifold in the parameter space. Nearby bifurcations should be studied using generic two-parameter families transverse to this codim 2 sub-manifold. We can therefore restrict ourselves to generic two-parameter ODEs ($m=2$), where the primary homoclinic orbit exists along a smooth \textit{homoclinic curve} in the parameter plane, while the saddle to saddle-focus transition happens at an isolated point on this curve. There are many more codim 2 homoclinic bifurcations, see \cite{ChaKuz1994,Shilnikovs2001,HomburgSandstede2010}.\\
 
As already noted in \cite{ChaKuz1994}, the simplest saddle to saddle-focus transitions correspond to
\begin{itemize}
\item[(i)] a \textit{double} leading eigenvalue;
\item[(ii)] \textit{three simple} leading eigenvalues.
\end{itemize}

\begin{figure}[htbp]
	\centering
	\begin{subfigure}[b]{0.25\textwidth}
		\includegraphics[width=\textwidth]{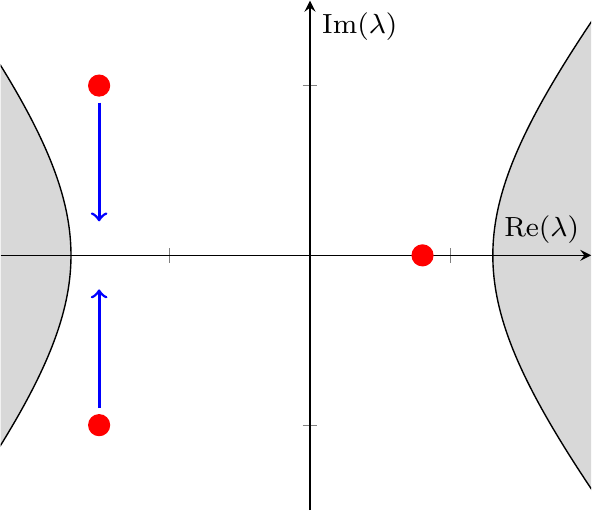}
		\caption{$\alpha<0$}
	\end{subfigure}\hfill
	\begin{subfigure}[b]{0.25\textwidth}
		\includegraphics[width=\textwidth]{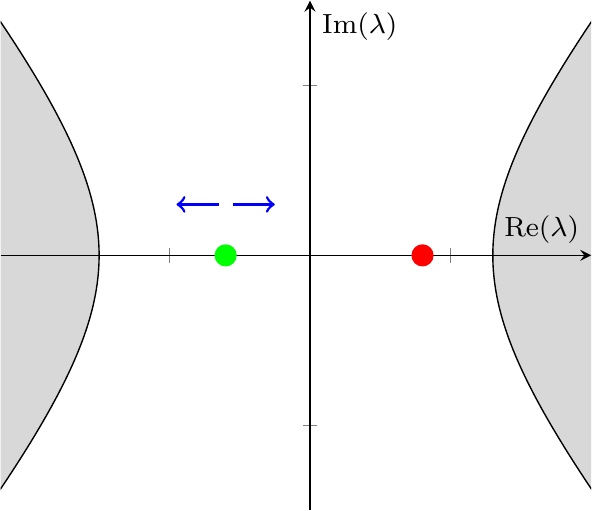}
		\caption{$\alpha=0$}
	\end{subfigure}\hfill
	\begin{subfigure}[b]{0.25\textwidth}
		\includegraphics[width=\textwidth]{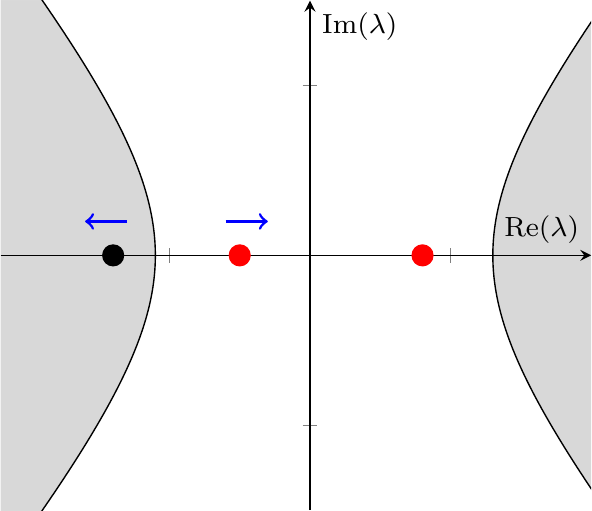}
		\caption{$\alpha>0$}
	\end{subfigure}
	\caption{Eigenvalue configurations of the saddle to saddle-focus transition in case (i); $\alpha$ is the parameter along the homoclinic curve and the bifurcation occurs at $\alpha=0$. Arrows point in the direction of generic movement of eigenvalues. The green marker indicates a double real eigenvalue. The gray areas  contain non-leading eigenvalues, leading eigenvalues are marked red and non-leading eigenvalues are marked black.}
	\label{fig:introbel}
\end{figure}
\begin{figure}[htbp]
	\centering
	\begin{subfigure}[b]{0.25\textwidth}
		\includegraphics[width=\textwidth]{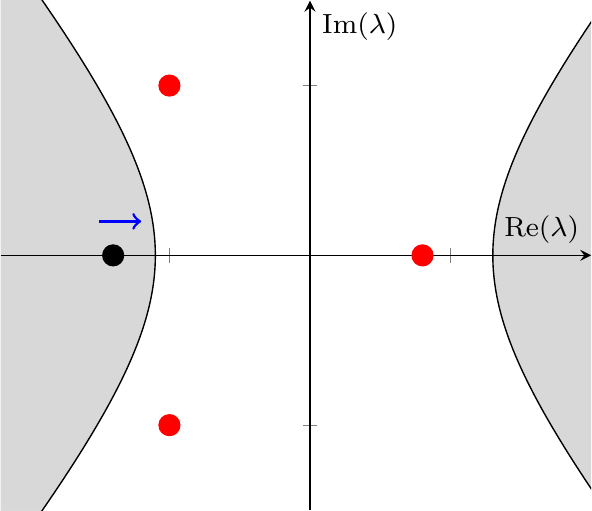}
		\caption{$\alpha<0$}
	\end{subfigure}\hfill
	\begin{subfigure}[b]{0.25\textwidth}
		\includegraphics[width=\textwidth]{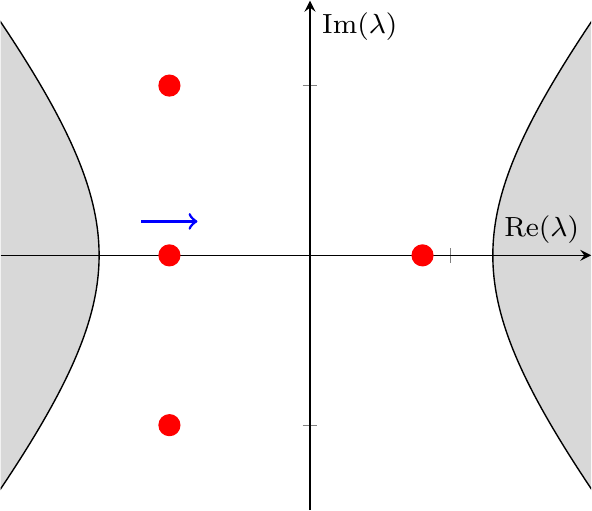}
		\caption{$\alpha=0$}
	\end{subfigure}\hfill
	\begin{subfigure}[b]{0.25\textwidth}
		\includegraphics[width=\textwidth]{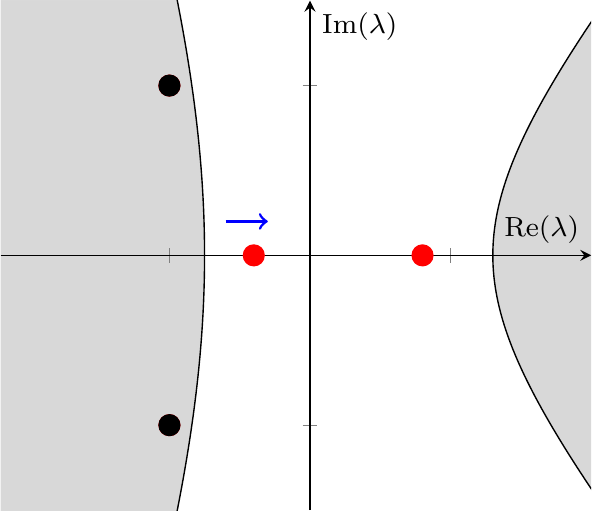}
		\caption{$\alpha>0$}
	\end{subfigure}
	\caption{Eigenvalue configurations of the saddle to saddle-focus transition in case (ii); the scalar bifurcation parameter along the homoclinic curve is $\alpha$. Arrows point in the direction of generic movement of eigenvalues. There is a codimension 2 situation at $\alpha=0$, where the leading stable eigenspace becomes 3-dimensional. Non-leading eigenvalues are contained in the gray area, leading eigenvalues are marked red and non-leading eigenvalues are marked black.}
	\label{fig:intro3DL}
\end{figure}

In case (i), see \autoref{fig:introbel}, the pair of leading complex eigenvalues approaches the real axis and splits into two distinct real eigenvalues. At the transition there is a double real eigenvalue and the leading eigenspace is two-dimensional. In case (ii), see \autoref{fig:intro3DL}, the real eigenvalue exchanges its position with the pair of complex-conjugate eigenvalues. At the transition there are two complex-conjugate eigenvalues and one real eigenvalue with the same real part. All leading eigenvalues are simple and the leading eigenspace is 3-dimensional.\\

Case (i) is a saddle to saddle-focus homoclinic transition that appears in various applications, e.g. in biophysics \cite{Kuznetsov1994} and ecology \cite{Kuznetsov2001a}. Moreover, in these applications the transition corresponds to the wild case with $\sigma_0>0$. 
This case was first studied theoretically by Belyakov \cite{Belyakov1980}, who proved that the corresponding bifurcation diagram is complicated. We call this case the {\em standard Belyakov case}. In \cite{Belyakov1980,Kuznetsov2001a} a description of the main features of the universal bifurcation diagram close to this transition for $n=3$ in the wild case has been obtained: 
\begin{enumerate}
	\item[1.] There exists an infinite set of \textit{limit point} (LP) and \textit{period doubling} (PD) bifurcation curves.
	\item[2.] There exists an infinite set of \textit{secondary homoclinic} curves corresponding to homoclinic orbits making two global excursions and various numbers of local turns near the equilibrium.  
   \item[3.] Both sets have the same `bunch' shape: The corresponding curves emanate from the codim 2 point and accumulate onto the branch of primary saddle-focus homoclinic orbits. The secondary homoclinics accumulate only from one side.
\end{enumerate}

Case (ii) has been observed in \cite{Meijer2014} for a 4D system of ODEs arising from a study of travelling waves in a neural field model. We note that the transition in \cite{Meijer2014} is tame with $\sigma_0<0$. As in the standard Belyakov case, we expect a complicated bifurcation diagram in the wild case, i.e. when $\sigma_0>0$. \\

Our paper is devoted to the theoretical analysis of the homoclinic saddle to saddle-focus transition for case (ii), when the leading stable eigenspace is three-dimensional, focusing on the wild case. We call this transition the \textit{3DL-transition}. To our best knowledge, no systematic analysis of this case is available in the literature. A possible reason for this gap is that case (ii) could only occur in (\ref{eq:intro1}) with $n \geq 4$, while case (i) happens already in three-dimensional ODEs. This leads to the study of a three-dimensional return map in case (ii), which is much more difficult to analyze than the planar return map in the standard Belyakov case (i).\\

For $n$-dimensional systems, generically the analysis of homoclinic bifurcations can be restricted to the \textit{Homoclinic Center Manifold}, a $k$-dimensional invariant finitely-smooth manifold that is tangent at the equilibrium to the eigenspace corresponding to the union of all \textit{leading} eigenvalues \cite{Homburg1996,Shilnikovs1998,ShaTur1999,Sandstede2000}. Thus, $k$ is the number of all \textit{leading} eigenvalues of the equilibrium, counting their multiplicities. For the saddle to saddle-focus homoclinic transition case (ii), we have $k=4$. Thus, the analysis of four-dimensional ODEs is sufficient to describe the main features of the bifurcation diagram near this transition in a generic $n$-dimensional situation.\\

By considering a generic 4D system with the 3DL-transition, we are able to obtain a two-parameter model 3D return map which describes the bifurcations occurring close to the transition. We will see that when $\sigma_0>0$, there exist infinitely many bifurcation curves. However, the shape of these bifurcation curves differs essentially from those in the standard Belyakov case (i):
\begin{enumerate}
	\item[1.] There exist infinitely many PD, LP, NS and secondary homoclinic curves. These curves accumulate onto the curve of primary homoclinic orbits but do not emanate from the codim 2 point.
	\item[2.] Each LP curve is a `horn' composed of two branches. Close to the horn's tipping point LP and PD curves interact via {\em spring} and {\em saddle areas} \cite{crossroad}. Transitions between saddle and spring areas are observed. Each secondary homoclinic curve forms a `horizontal parabola'. 
	\item[3.] Several codim 2 points exist on each of the LP, PD, and NS curves. We observe GPD (along PD) and cusp (along LP) points, as well as strong resonances. 
\end{enumerate}
    
Using the model map, we prove analytically that the cusp points asymptotically approach the 3DL transition point. The same is shown for the secondary homoclinic turning points. We present numerical evidence that all other mentioned codimension 2 points form sequences also converging to the 3DL transition point. \\

This paper is organized as follows. In Section \ref{sec:return_map} we formulate the genericity assumptions on (\ref{eq:intro1}) with $n=4$ and $m=2$. Next,  we derive a model 3D return map and its 1D simplification. The $C^1$-linearization theorem by Belitskii is used near the equilibrium. In Section \ref{sec:scalar_map} we analyze the 1D model map to describe LP and PD bifurcations of the fixed points/periodic orbits. An essential part of the analysis of the 1D map is carried out analytically, while that of the full 3D model map in Section \ref{sec:3D_map} employs advanced numerical continuation tools, except for the LP and PD bifurcations (reducible to the 1D return map studied in Section \ref{sec:scalar_map}) and the secondary homoclinic bifurcations. In Section \ref{sec:ODE_interpretation}, implications for the dynamics of the original 4D ODE system  are summarized.  Finally, in Section \ref{sec:example_model}, we give an explicit example of the chaotic 3DL transition in a perturbed Lorenz-Stenflo model.  

\section{Derivation of the model maps}
\label{sec:return_map}

\subsection{Assumptions}
We make the following assumptions about the 3DL-transition at the critical parameter values, which we assume to be $\alpha_1=\alpha_2=0$. Recall that we only consider $n=4$ and $m=2$.

\begin{enumerate}
	\item[] \textbf{(A.1)} The eigenvalues of the linearisation at the critical 3DL equilibrium $x=0$ are  
    
$$
\delta_0,\delta_0 \pm i \omega_0 \textnormal{ and } \epsilon_0,
$$  
where $\delta_0<0, \omega_0>0$ and $\epsilon_0>0$.
	\item[] \textbf{(A.2)} There exists a primary homoclinic orbit $\Gamma_0$ to this 3DL equilibrium. 
	\item[] \textbf{(A.3)} The homoclinic orbit $\Gamma_0$ satisfies the following genericity condition: The normalized tangent vector to $\Gamma_0$ has nonzero projections to both the 1D eigenspace corresponding to the real eigenvalue $\delta_0$ and to the 2D eigenspace corresponding to the complex eigenvalues $\delta_0 \pm i \omega_0$, when approaching the equilibrium.  	
\end{enumerate}

 
Any system satisfying the assumptions \textbf{(A.1-3)}, can be transformed near the critical equilibrium via a translation, a linear transformation, a linear time scaling, and introducing new parameters $\mu=(\mu_1,\mu_2)$, to 
\begin{equation}
\left\{ \begin{array}{r@{\quad}cr}
\dot{x_1} = & \gamma(\mu) x_1 - x_2 + \tilde{f}_1(x,\mu), \\ 
\dot{x_2} = &  x_1 + {\gamma}(\mu) x_2 + \tilde{f}_2(x,\mu),\\ 
\dot{x_3} = &\left({\gamma}(\mu) - \mu_1 \right) x_3 +\tilde{f}_3(x,\mu),\\ 
\dot{x_4} = & {\beta}(\mu)x_4 + \tilde{f}_4(x,\mu), 
\end{array}\right.
\label{eq:system}
\end{equation}
where the functions $\tilde{f}_i$ are nonlinear and such that $\tilde{f}_i(0,\mu) = 0$ for $i=1,2,3,4$ for all $\mu \in\R^2$ sufficiently small,
\begin{equation}
{\gamma}(0) = \frac{\delta_0}{\omega_0} \textnormal{\ \ and\ \ } {\beta}(0) = \frac{\epsilon_0}{\omega_0}. 
\label{eq:eigsat3DL}	
\end{equation} 

Here $\mu_2=\mu_2(\alpha)$ is a `splitting function' so that the primary homoclinic orbit to the equilibrium (saddle, 3DL, saddle-focus) exists along the curve $\mu_2=0$. The exact choice of $\mu_2$ will be clear later. The value $\mu_1(\alpha)$ controls which stable eigenvalue leads. For $\mu_1>0$, the stable leading eigenvalues are complex (saddle-focus case) and for $\mu_1<0$ the stable leading eigenvalue is real (saddle case).\\

Now we can formulate the final (transversality) assumption:
\begin{enumerate}
	\item[] \textbf{(A.4)}  The components of $\mu=(\mu_1,\mu_2)$ are small and the 3DL saddle exists at $\mu=0$. Moreover, the mapping $\alpha \mapsto \mu(\alpha)$ is regular at $\alpha=0$.
\end{enumerate}    
 
\subsection{$C^1$-linearisation near the equilibrium}
For  $\mu_1$ sufficiently small, we can use Belitskii's Theorem \cite{Belitskii1979,DS-5} to get the $C^1$-equivalence of the flow generated by \eref{eq:system} to that corresponding to its linear part, near the equilibrium $O=(0,0,0,0)$. Therefore, near $O$, we can use the following linear system:
\begin{equation}
\left\{ \begin{array}{r@{\quad}cr}
\dot{x_1} = & \gamma(\mu) x_1 - x_2, \\ 
\dot{x_2} = &  x_1 + \gamma(\mu) x_2,\\ 
\dot{x_3} = &\left(\gamma(\mu) - \mu_1 \right) x_3,\\ 
\dot{x_4} = & \beta(\mu)x_4.\\ 
\end{array}\right.
\label{eq:systemlin}
\end{equation}
Outside of a fixed neighborhood of $O$, the nonlinear terms should be taken into account.

\subsection{Introducing cross-sections}
Our next aim is to derive the model Poincar\'e map close to $\Gamma_0$ near the 3DL transition, that we will use for the two-parameter perturbation study.\\

\begin{figure}[hbtp]
	\centering
	\includegraphics[width=0.45\textwidth]{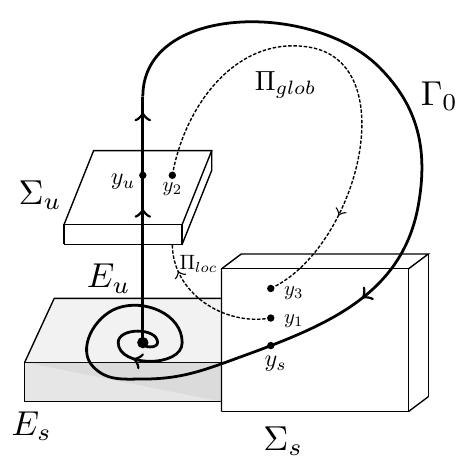}
	\caption[Geometric construction of the near-to-saddle map]{The geometric construction of cross-sections close to the critical 3DL-saddle at $(0,0,0,0)$ and the homoclinic connection $\Gamma_0$, in order to obtain the map $\Pi:\Sigma_s \mapsto \Sigma_s$. Here $\Sigma_u$ is defined by the cross-section $x_4=1$ and $\Sigma_s$ is the cross-section $x_3=0$. The homoclinic connection is assumed to pass through the points $y_s=(1,0,1,0)$ and $y_u = (0,0,0,1)$. The stable and unstable eigenspaces are $E_s$ and $E_u$, respectively. }
	\label{fig:poincare_map}
\end{figure}

\autoref{fig:poincare_map} gives an impression of the homoclinic connection to a 3DL-saddle in the four-dimensional system. As we are interested in understanding the bifurcations close to the saddle and the homoclinic orbit, we define two Poincar\'e cross-sections,
\begin{eqnarray}
\Sigma_s &= \left\{ (x_1,x_2,x_3,x_4) | x_2=0 \right\},\\
\Sigma_u &= \left\{ (x_1,x_2,x_3,x_4) | x_4=1 \right\},
\end{eqnarray}  

and assume, after linearly rescaling variables if necessary, that the homoclinic orbit passes through these cross-sections at $y_s = (1,0,1,0)$ and $y_u = (0,0,0,1)$ respectively, for all parameter values along the primary homoclinic curve. \\

Clearly, both cross-sections are transversal to the flow and to the stable and unstable eigenspaces. Thus, by following orbits starting from $\Sigma_s$ to $\Sigma_u$ and returning back to $\Sigma_s$, we will be able to define a three-dimensional map $\Pi$ mapping (a subset of) $\Sigma_s$ to itself. This map can then be used to study both periodic orbits and secondary homoclinic orbits.\\

We shall construct the map $\Pi$ by composing two maps, $\Pi_{loc}:\Sigma_s \mapsto\Sigma_u$ and $\Pi_{glob}:\Sigma_u \mapsto \Sigma_s$, i.e.
\begin{eqnarray}
\Pi = \Pi_{glob} \circ \Pi_{loc}.
\end{eqnarray}
\subsection{Derivation of the return map}
\label{sec:Derivation}
We begin by computing the orbit of \eref{eq:systemlin} starting from an arbitrary point  in $\Sigma_s$, close to $y_s$. This will be used to define the map $\Pi_{loc}$. The local map $\Pi_{loc}$ is the mapping from  $y_1=(x_1^s,0,x_3^s,x_4^s) \in \Sigma_s$ to a point $y_2=(x_1^u,x_2^u,x_3^u,1) \in \Sigma_u$ close to $y_u$, and is given by

\begingroup
\renewcommand*{\arraystretch}{1.8}
\begin{equation}
\Pi_{loc}:
\left(\begin{array}{l}
x_1^s \\ 
x_3^s \\ 
x_4^s
\end{array}\right) \mapsto 	
\left(\begin{array}{l}
x_1^s (x_4^s)^{\nu(\mu)} \cos\left(-\frac{1}{\beta(\mu)}\ln{x_4^s}\right)\\
x_1^s (x_4^s)^{\nu(\mu)} \sin\left(-\frac{1}{\beta(\mu)}\ln{x_4^s}\right)\\
x_3^s(x_4^s)^{\nu(\mu)+\mu_1/\beta(\mu)} 
\end{array}\right).	
\label{eq:locmap}
\end{equation}
\endgroup 
where $\nu(\mu) = -\gamma(\mu)/\beta(\mu)$. Let $\gamma_0:=\gamma(0)$ and $\beta_0:=\beta(0)$. The quantity 
\begin{equation}
\label{saddle_index}
\nu_0=-\frac{\gamma_0}{\beta_0}=-\frac{\delta_0}{\epsilon_0}
\end{equation}
is called the \textit{saddle index}. Note that the saddle quantity $\sigma_0$ is related to the saddle index \eref{saddle_index} as follows:
\begin{eqnarray*}
\nu_0<1 \Longleftrightarrow \sigma_0>0,\\
\nu_0>1 \Longleftrightarrow \sigma_0<0.
\end{eqnarray*}
We assume from now on that $\nu_0<1$, so that only the wild case $\sigma_0>0$ is considered.\\

For the global return map $\Pi_{glob}:\Sigma_u \mapsto \Sigma_s$, we obviously cannot use \eref{eq:systemlin}. Instead  we use a general $C^1$ approximation of the flow from $(0,0,0,1)$ to $(1,0,1,\mu_2)$. Here $\mu_2$ is the aforementioned splitting parameter. It controls the return of the orbit to the critical saddle. For $\mu_2=0$ only, we have a primary homoclinic connection.\\

Thus, the following representation of $\Pi_{glob}$ can be used
\begingroup
\renewcommand*{\arraystretch}{1.8}
\begin{equation}
\hspace{-1.5cm}\Pi_{glob}:
\left(\begin{array}{l}
x_1^u \\ 
x_2^u \\ 
x_3^u
\end{array}\right) \mapsto	
\left(\begin{array}{l}
1 \\
1 \\
\mu_2
\end{array}\right)	+
\left(\begin{array}{lll}
a_{11}(\mu) & a_{12}(\mu) & a_{13}(\mu) \\
a_{21}(\mu) & a_{22}(\mu) & a_{23}(\mu) \\
a_{31}(\mu) & a_{32}(\mu) & a_{33}(\mu)
\end{array}\right) 
\left(\begin{array}{l}
x_1^u \\
x_2^u \\
x_3^u
\end{array}\right)+
o(\|x^u\|),
\label{eq:globmap}
\end{equation}
\endgroup

where $x^u=(x_1^u,x_2^u,x_3^u)$. For $A_0 = [a_{ij}(0)]$ we also have  $\det(A_0) \neq 0$ which follows from the invertibility of $\Pi_{glob}$ for $\mu$ small enough.\\

Now, \eref{eq:locmap} and \eref{eq:globmap} together give us the full return map $\Pi = \Pi_{glob} \circ \Pi_{loc}$. Keeping the dependence of all coefficients on $\mu$ implicit, we can write $\Pi$ as
\begingroup
\renewcommand*{\arraystretch}{1.8}
\begin{equation}
\hspace{-2.5cm}\Pi:\left(\begin{array}{l}
x_1^s \\ 
x_3^s \\ 
x_4^s
\end{array}\right) \mapsto 	
\left(\begin{array}{l}
1 + b_1 x_1^s (x_4^s)^\nu \cos \left( -\frac{1}{\beta}\ln x_4^s + \theta_1 \right) + b_2 x_3^s (x_4^s)^{\nu+\mu_1/\beta}\\
1 + b_3 x_1^s (x_4^s)^\nu \sin \left( -\frac{1}{\beta}\ln x_4^s + \theta_2 \right) + b_4 x_3^s (x_4^s)^{\nu+\mu_1/\beta} \\
\mu_2 + b_5 x_1^s (x_4^s)^\nu \sin \left( -\frac{1}{\beta}\ln x_4^s + \theta_3 \right) + b_6 x_3^s (x_4^s)^{\nu+\mu_1/\beta}
\end{array}\right)+o(\|x^s\|^{\nu}),
\label{eq:finalmap1}
\end{equation}
\endgroup
where $x^s=(x_1^s,x_3^s,x_4^s)$ and
\begingroup
\renewcommand*{\arraystretch}{1.8}
\begin{equation}
\begin{array}{lll}
\sin \theta_1 = -\frac{a_{12}}{\sqrt{a_{11}^2+a_{12}^2}}, & \cos \theta_2 = \frac{a_{22}}{\sqrt{a_{21}^2+a_{22}^2}}, & \cos \theta_3 = \frac{a_{32}}{\sqrt{a_{31}^2+a_{32}^2}}, \\
b_1 = \sqrt{a_{11}^2 + a_{12}^2}, & b_3 = \sqrt{a_{21}^2 + a_{22}^2}, & b_5 = \sqrt{a_{31}^2 + a_{32}^2},\\
b_2 = a_{13}, & b_4=a_{23}, & \textnormal{and } \ b_6 = a_{33}.
\end{array}
\end{equation}
\endgroup

Following \ref{Gonchenko1997}, we make the smooth invertible transformation $x_4^s \mapsto x_4^s\exp \left( \theta_3 \beta \right)$ to eliminate $\theta_3$. This gives
\begingroup
\renewcommand*{\arraystretch}{1.8}
\begin{equation}
\hspace{-2.0cm}\Pi:\left(\begin{array}{l}
x_1 \\ 
x_3 \\ 
x_4
\end{array}\right) \mapsto 	
\left(\begin{array}{l}
1 + \alpha_1 x_1 x_4^\nu \cos \left( -\frac{1}{\beta}\ln x_4 + \phi_1 \right) + \alpha_2 x_3 x_4^{\nu+\mu_1/\beta}\\
1 + \alpha_3 x_1 x_4^\nu \sin \left( -\frac{1}{\beta}\ln x_4 + \phi_2  \right) + \alpha_4 x_3 x_4^{\nu+\mu_1/\beta} \\
\mu_2 + C_1 x_1 x_4^\nu \sin \left( -\frac{1}{\beta}\ln x_4  \right) + C_2 x_3 x_4^{\nu+\mu_1/\beta}
\end{array}\right)+o(\|x\|^{\nu}),
\label{eq:3DmapFull}
\end{equation}
\endgroup
where we have dropped the superscript `$s$' from the coordinates of $x=(x_1,x_3,x_4)$ for convenience,  and where
\begin{eqnarray}
\begin{array}{ll}
\phi_1 = \theta_1-\theta_3, & \phi_2 = \theta_2-\theta_3, \\ \alpha_1 = b_1 \exp (\theta_3 \beta \nu), &
 \alpha_2 = b_2 \exp((\nu+\mu_1/\beta)\theta_3\beta) , \\ \alpha_3 = b_3 \exp (\theta_3 \beta \nu), & \alpha_4 =  b_4 \exp((\nu+\mu_1/\beta)\theta_3\beta) , \\ C_1=b_5 \exp (\theta_3 \beta \nu) & \textnormal{and } \ C_2 = b_2 \exp((\nu+\mu_1/\beta)\theta_3\beta). 
\end{array}
\end{eqnarray}
Observe that all constants $\alpha_j$ and $C_k$ depend on $\mu$ and that $C_1>0$. Let us denote by $\alpha_j^0$ and $C_j^0$ their critical values at $\mu=0$.\\

Truncating the $o(\|x\|^{\nu})$-terms in (\ref{eq:3DmapFull}) and taking only the critical values of all constants, we define 
\begin{equation}
G(x,\mu):=	
\left(\begin{array}{l}
1 + \alpha^0_1 x_1 x_4^{\nu_0} \cos \left( -\frac{1}{\beta_0}\ln x_4 + \phi^0_1 \right) + \alpha^0_2 x_3 x_4^{\nu_0+\mu_1/{\beta_0}}\\
1 + \alpha^0_3 x_1 x_4^{\nu_0} \sin \left( -\frac{1}{\beta_0}\ln x_4 + \phi^0_2  \right) + \alpha^0_4 x_3 x_4^{\nu_0+\mu_1/{\beta_0}} \\
\mu_2 + C^0_1 x_1 x_4^{\nu_0} \sin \left( -\frac{1}{\beta_0}\ln x_4  \right) + C^0_2 x_3 x_4^{\nu_0+\mu_1/{\beta_0}}
\end{array}\right).
\label{eq:3Dmap}
\end{equation}
This map $G$ is the final form of the \textit{3D model return map} that we will use for the numerical analysis ahead.\\

Now, to analyze periodic orbits close to the homoclinic connection with respect to the critical 3DL-saddle, we look for fixed points of the map $\Pi$ given by \eref{eq:3DmapFull}. These fixed points correspond to periodic orbits in the original ODE system. Bifurcations of these fixed points describe the various local bifurcations of the corresponding periodic orbits. \\

The fixed point condition for map \eref{eq:3DmapFull} is
\begingroup
\renewcommand*{\arraystretch}{1.8}
\begin{equation}
\hspace{-1.0cm}\left(\begin{array}{l}
x_1 \\ 
x_3 \\ 
x_4
\end{array}\right) = 	
\left(\begin{array}{l}
1 + \alpha_1 x_1 x_4^{\nu} \cos \left(-\frac{1}{\beta}\ln x_4 + \phi_1 \right) + \alpha_2 x_3 x_4^{\nu+\mu_1/{\beta}}\\
1 + \alpha_3 x_1 x_4^{\nu} \sin \left(-\frac{1}{\beta}\ln x_4 + \phi_2  \right) + \alpha_4 x_3 x_4^{\nu+\mu_1/{\beta}} \\
\mu_2 + C_1 x_1 x_4^{\nu} \sin \left(-\frac{1}{\beta}\ln x_4  \right) + C_2 x_3 x_4^{\nu+\mu_1/{\beta}}
\end{array}\right)+o(\|x\|
^{\nu}),
\label{eq:fixedpt1}
\end{equation}
\endgroup
where all constants $\alpha_j$ and $C_k$ still depend on $\mu$. For non-degeneracy, we require that real constants $C_1$ and $C_2$ are nonzero. We justify this later. Substituting the expressions for $x_1$ and $x_3$ from \eref{eq:fixedpt1} in the expression for $x_4$, we get
\begin{eqnarray}
x_4 = \mu_2 + C_1 x_4^\nu \sin \left( -\frac{1}{\beta}\ln x_4  \right) + C_2  x_4^{\nu+\mu_1/\beta} + o(|x_4|^{\nu}),
\label{eq:1Dmap1}
\end{eqnarray}

as our one-dimensional fixed point condition. As we are interested behavior close to $(1,0,1,0)$ on the cross-section $\Sigma_s$, we consider only the leading terms of \eref{eq:1Dmap1} and introduce the following \textit{scalar model map}:
\begin{eqnarray}
F(x,\mu):=  \mu_2 + C^0_1 x^{\nu_0} \sin \left(- \frac{1}{\beta_0}\ln x  \right) + C^0_2  x^{\nu_0+\mu_1/{\beta_0}}.
\label{eq:1dmap0}
\end{eqnarray}
The extra additive term $C^0_2  x^{\nu_0+\mu_1/{\beta_0}}$ is what makes this map different from the scalar model maps describing the codim 1 saddle-focus case.\\

If we were to set $C^0_1$ to zero, then we would obtain finitely many fixed points for all values of $\nu_0, \mu_1, \beta_0, \mu_2$ and $C^0_2$. If we set $C^0_2$ to zero, we get the codim 1 saddle-focus case. Thus we assume
\begin{enumerate}
	\item[] \textbf{(A.5)}  $C^0_1C^0_2 \neq 0.$
\end{enumerate}    

\section{Analysis of the scalar model map}
\label{sec:scalar_map}
In this section, we study bifurcations of fixed points of the map \eref{eq:1dmap}. To stay close to the 3DL bifurcation, we only work with small values of $x$ and $\mu$. To simplify notations, we rewrite the scalar model map as
\begin{eqnarray}
F(x,\mu):=  \mu_2 + C_1 x^{\nu} \sin \left(- \frac{1}{\beta}\ln x  \right) + C_2  x^{\nu+\mu_1/{\beta}},
\label{eq:1dmap}
\end{eqnarray}
assuming that $\nu, \beta$, and $C_{1,2}$ are fixed at their critical values.

\begin{figure}[b!]
	\centering
	\begin{subfigure}[b]{0.45\textwidth}
		\includegraphics[width=\textwidth]{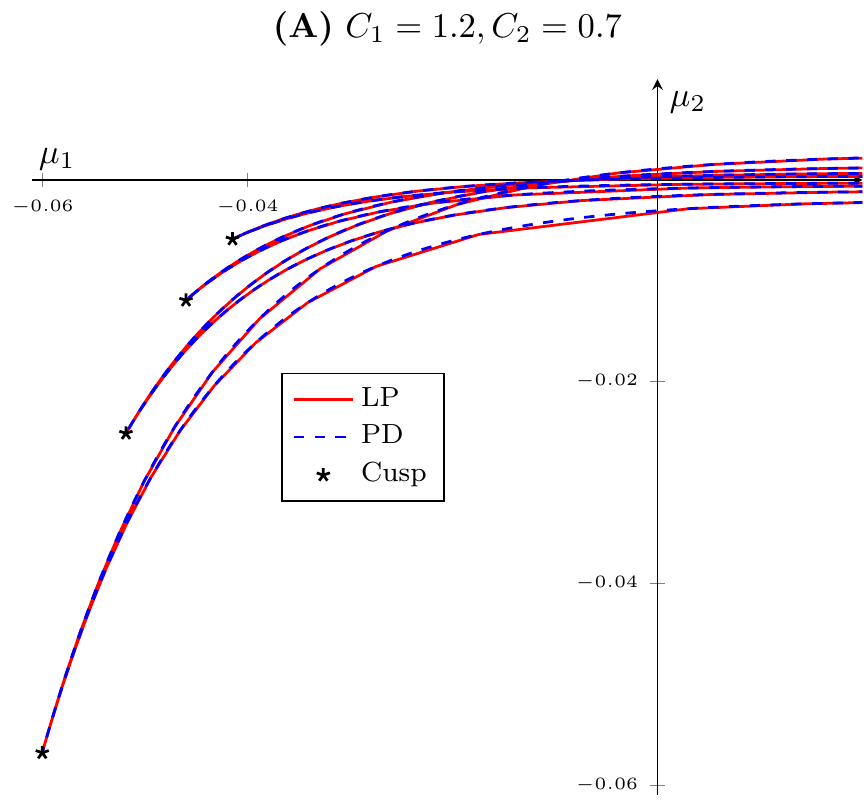}
	\end{subfigure}\hfill
	\begin{subfigure}[b]{0.43\textwidth}
		\includegraphics[width=\textwidth]{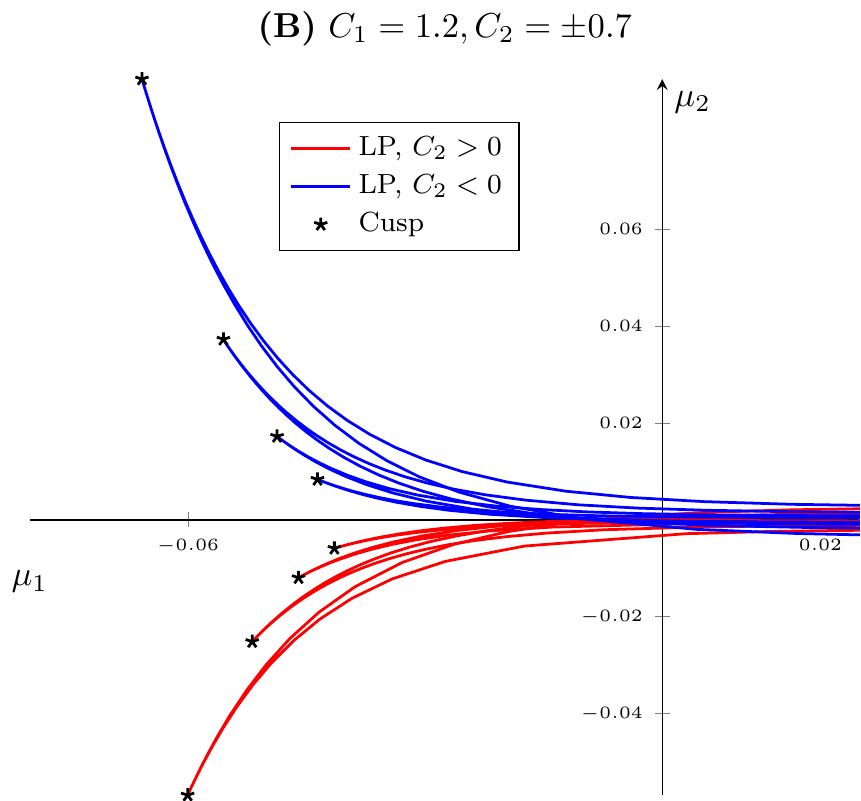}
	\end{subfigure}\\
	\begin{subfigure}[b]{0.45\textwidth}
		\includegraphics[width=\textwidth]{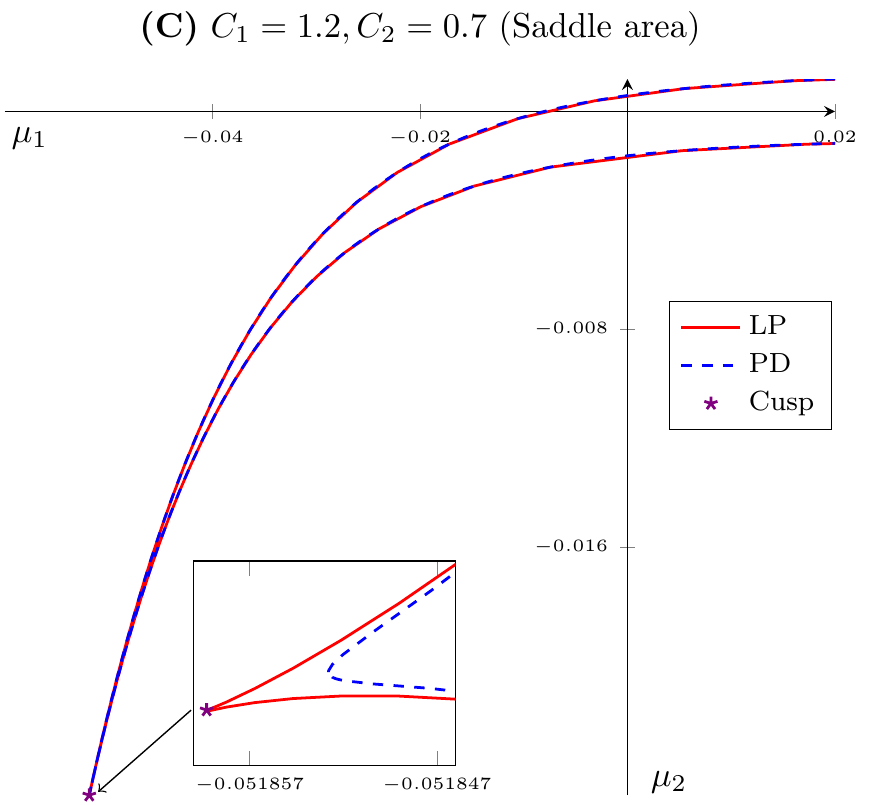}
	\end{subfigure}	\hfill
	\begin{subfigure}[b]{0.45\textwidth}
		\includegraphics[width=\textwidth]{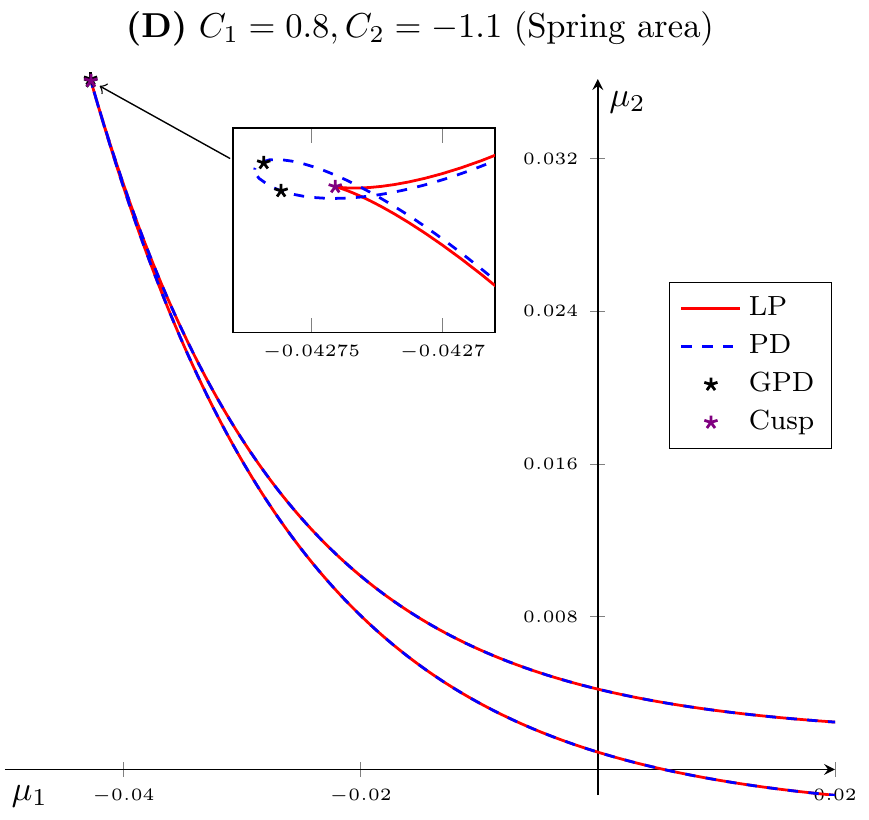}
	\end{subfigure}
	\caption{Primary LP and PD bifurcation curves obtained by numerical continuation, for the map \eref{eq:1dmap} for some representative values of $C_1$ and $C_2$. In panel \textbf{(A)} we plot 4 pairs of these curves. All of them have the same global structure. There are two types of codimension 2 points that can be found along these curves: Cusp (on LP curves) and GPD (along PD curves). In panel \textbf{(B)} we see what happens when we switch the sign of $C_2$, the \textit{horns} move from $\mu_2>0$ to $\mu_2<0$. In panel \textbf{(C)} and \textbf{(D)} we see examples of one PD and LP curve with the saddle area and spring area respectively (zoomed in). In the insets, $\mu_2$ is scaled for convenience.}
	\label{fig:1Dcont}
\end{figure}
\subsection{Numerical continuation results}
Using the continuation package MatcontM \cite{matcontm1,matcontm2}, we obtained many LP and PD bifurcation curves, which form interesting structures. There is strong evidence that there exist infinitely many PD and LP curves in the $(\mu_1,\mu_2)$-parameter space. Several such curves can be seen in \autoref{fig:1Dcont}. We make the following observations:
\begin{enumerate}
	\item   The curves exhibit a repetitive behavior: two branches of one LP curve meet to form a \textit{horn}. Apparently, infinitely many such `horns' exist. The sequence of these `horns' in the parameter space appears to approach $\mu_2=0$ asymptotically, which is the curve of primary homoclinic orbits. Also, the tips of the `horns' are always located entirely in either the second, or third quadrant of the $(\mu_1,\mu_2)$-space.
	\item  The PD and LP curves appear to coincide on visual inspection, and there can exist GPD points in the vicinity of the tip of the LP horn. 
	\item  The tip each LP `horn' is a cusp point. These cusps always exist, for all values of $C_1$ and $C_2$ and form a sequence that appears to approach the origin $\mu=0$.
	\item  Upon closer inspection, we observe that there exists either of the two subtle structures near the top of every LP `horn'. One is a \textit{spring area}, where the PD curve loops around the cusp point. The other is a \textit{saddle area}, where the PD curve makes a sharp turn close to the cusp, see the insets in \autoref{fig:1Dcont}. The spring area is accompanied by two GPD points along the PD loop. These points are absent in a saddle area. Mira et al. \cite{crossroad} discuss in detail the spring and saddle areas, including transitions from one case to the other and their genericity.    
\item The global behavior of this set of curves depends on parameters $C_1$ and $C_2$. For example, by switching the sign of $C_2$, the set of curves can be moved from the second to the third quadrant of the $\mu$-space, or vice-versa. The presence of saddle or spring areas depend on the parameters $C_1$ and $C_2$ but the exact conditions are not clear.
\end{enumerate}

In the sections ahead we support most of the observations by looking at analytical asymptotics of the LP and PD bifurcation curves of \eref{eq:1dmap}.
\begin{figure}[hbtp]
	\centering
	\begin{subfigure}[b]{0.6\textwidth}
		\includegraphics[width=\textwidth]{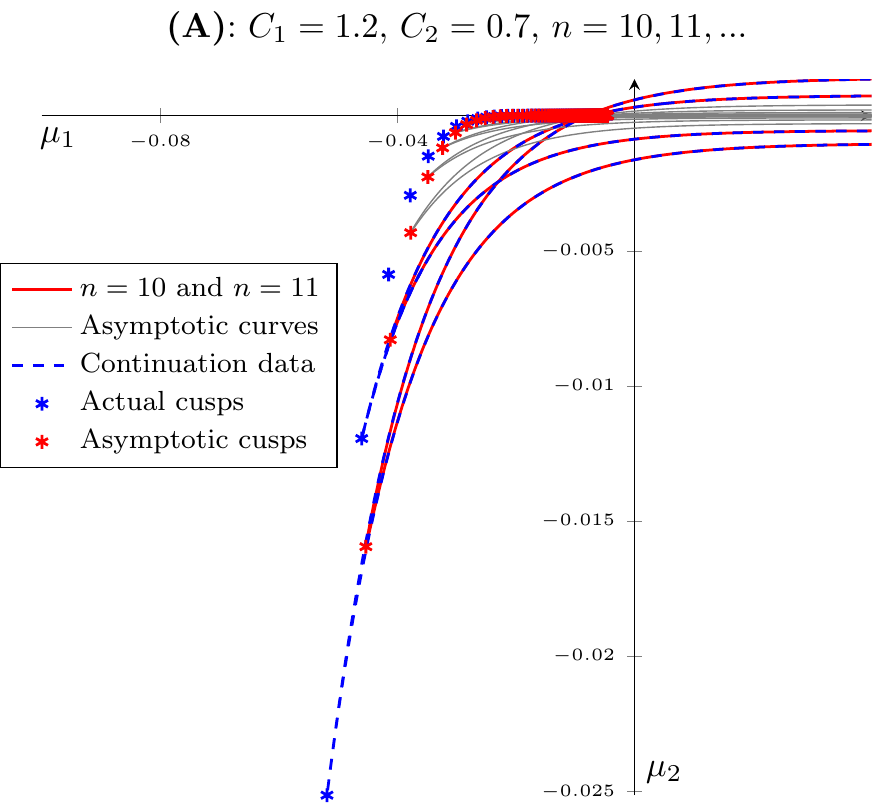}%
	\end{subfigure}\hfill
	\begin{subfigure}[b]{0.4\textwidth}
		\includegraphics[width=\textwidth]{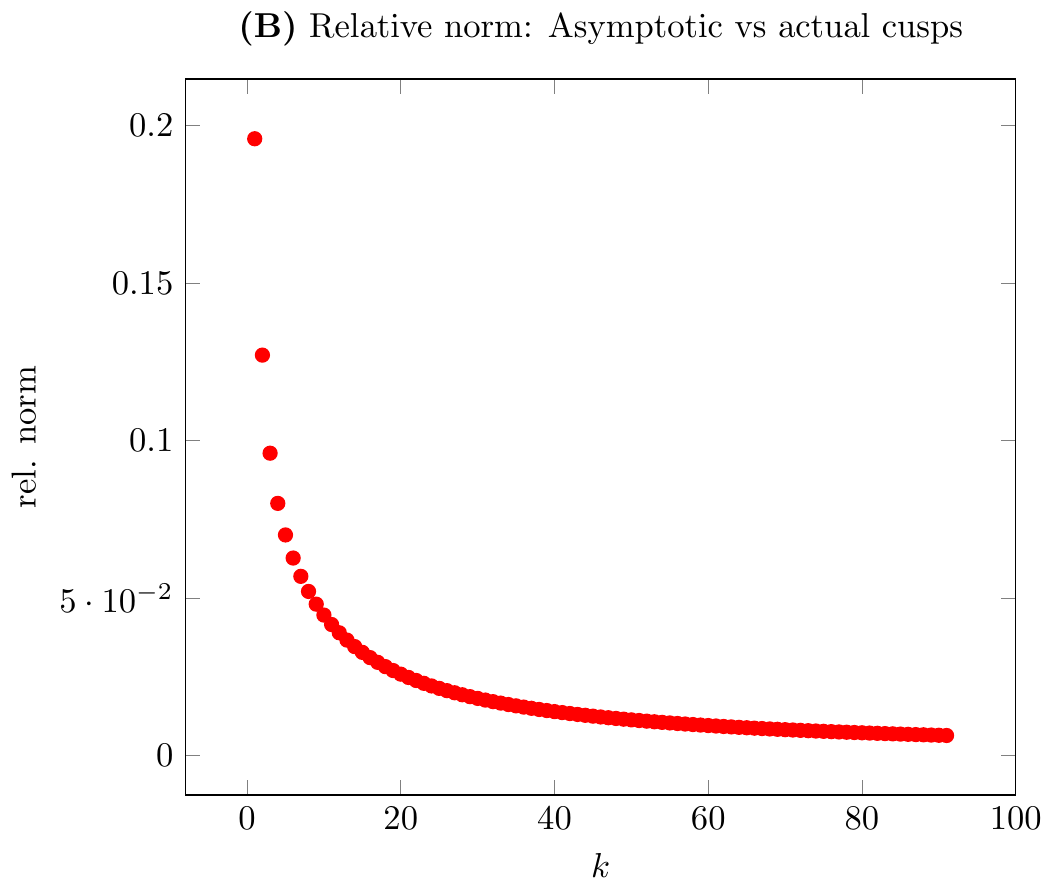}%
	\end{subfigure}
	\caption{Plots of asymptotic curves and PD/LP curves (obtained from continuation). In \textbf{(A)} we see how successive asymptotic curves in $n$ approximate the set of PD/LP curves. Here, cusps are obtained by performing Newton iterations to the defining system of the cusp bifurcation with starting points as the asymptotic cusps. In \textbf{(B)}, convergence of the asymptotic cusps to the actual cusps is observed. The corresponding values of $n$ in both plots are $n=10,11...,90$.  }
	\label{fig:AsympVsReal}
\end{figure}
\subsection{Asymptotics}
In this section we derive approximate solutions to the LP and PD conditions, and use them to justify numerical observations. As we are interested in solutions close to the 3DL bifurcation point $(\mu_1,\mu_2)=(0,0)$ we assume that $x,\mu_1$ and $\mu_2$ are sufficiently small. As we investigate only the wild case we restrict ourselves to $\nu<1$.
\subsubsection{LP horns and cusp points}
For the scalar model map \eref{eq:1dmap} the fixed point condition is given by
\begin{equation}
\mu_2 +C_1x^\nu \sin \left(-\frac{1}{\beta} \ln x \right) + C_2 x^{\nu+\mu_1/\beta} -x=0.
\label{eq:lpcond}
\end{equation} 
We parametrize $x$ using the following relation,
\begin{equation}
-\frac{1}{\beta} \ln x = 2 \pi n + \theta,
\end{equation}
for large $n \in \mathbb{N}$  and $\theta \in (0,2\pi)$. Thus, \eref{eq:lpcond} becomes
\begin{equation}
\mu_2 + C_1 e^{-\beta \nu (2\pi n +\theta)} \sin \theta + C_2 e^{-\beta (\nu+\mu_1/\beta)(2\pi n + \theta)} -e^{-\beta(2 \pi n+\theta)}=0.
\label{eq:lpcond2}
\end{equation} 
Let us define 
\begin{equation}
\hspace{-1.2cm}
	\Phi(\theta,\mu_1,\mu_2) := \mu_2 + C_1 e^{-\beta \nu (2\pi n +\theta)} \sin \theta + C_2 e^{-\beta (\nu+\mu_1/\beta)(2\pi n + \theta)} - e^{-\beta(2 \pi n+\theta)}.
\end{equation}
Then 
$$
\Phi_\theta(\theta,\mu_1,\mu_2) = 0,
$$
is the extra condition for an LP point. Computing the derivative, we get
\begin{equation}
\hspace{-1.2cm}
	C_1\left(\beta\nu \sin \theta - \cos \theta \right) + C_2(\beta\nu+\mu_1)e^{-\mu_1(2\pi n +\theta)} - \beta e^{-\beta(1-\nu)(2\pi n +\theta)}=0.
	\label{eq:lpcond3}
\end{equation}
We now simultaneously solve \eref{eq:lpcond} and \eref{eq:lpcond3} to obtain a sequence of functions $\mu_2^{(n)}(\mu_1)$ which describe the sequence of LP `horns' already observed numerically. Thus, rewriting \eref{eq:lpcond3}, we have  
\begin{equation}
\hspace{-1.2cm}
\begin{array}{ll}
 \beta \nu \sin \theta -\cos \theta  &= -\frac{C_2}{C_1}\left(\beta\nu+\mu_1\right)e^{-\mu_1(2\pi n +\theta)} +\frac{\beta}{C_1} e^{-\beta(1-\nu)(2\pi n +\theta)}\\
                                     &= -\frac{C_2}{C_1}\left(\beta \nu+\mu_1\right)e^{-2\pi\mu_1 n}\left[1-\mu_1\theta + O(\mu_1^2)\right] + O( e^{-\alpha n})\\	
                                     &= -\frac{C_2}{C_1}e^{-2\pi\mu_1 n}\left[\beta \nu - (1-\beta \nu \theta)\mu_1 + O(\mu_1^2)\right] + O( e^{-\alpha n}),
	\end{array}
\end{equation}
where $\alpha = 2 \pi \beta (1-\nu)>0$. Collecting trigonometric terms on the left we get
\begin{equation}
\hspace{-1.5cm}
	\sin(\theta-\phi)  =-\frac{1}{\sqrt{1+\beta^2\nu^2}}\frac{C_2}{C_1}e^{-2\pi\mu_1 n}\left[\beta \nu - (1-\beta \nu \theta)\mu_1 + O(\mu_1^2)\right] + O( e^{-\alpha n}),
	\label{eq:trigterms}
\end{equation}
where $\sin \phi = (1+\beta^2\nu^2)^{-1/2}$ and $\phi \in (0,\pi/2)$. Note that for large $n$ and negative $\mu_1$, the corresponding solution $\theta$ exists only for small $|\mu_1|$. Let
\begin{equation}
	\theta^n_0 := \arcsin \left( -\frac{\beta \nu}{\sqrt{1+\beta^2\nu^2}}\frac{C_2}{C_1}e^{-2\pi\mu_1 n}\right). 
	\label{eq:theta0}
\end{equation}
Then  we have two solutions,
\begin{equation}
\begin{array}{l}
\theta_1 = \phi + \theta^n_0 + \delta_{i1}(2\pi) + O(\mu_1),\\
\theta_2 = \phi + \pi - \theta^n_0 + O(\mu_1),
\label{eq:thetapdlp}
\end{array}
\end{equation} 
where $i = -\textnormal{sign}(C_2)$ and $\delta_{ij}$ is the Kronecker delta. \\

For each $n$ we obtain two solutions $\theta$ given by \eref{eq:thetapdlp}. The corresponding functions $\mu_2^{(n)}(\mu_1)$ are obtained from the first equation of \eref{eq:lpcond2},
\begin{equation}
\hspace{-1.5cm}
\left\{
\begin{array}{l}
	\mu_2^{(n,1)}(\mu_1) = -C_1 e^{-\beta \nu (2\pi n +\theta_1)} \sin \theta_1 - C_2 e^{-\beta (\nu+\mu_1/\beta)(2\pi n + \theta_1)} + e^{-\beta(2 \pi n+\theta_1)},\\
	\mu_2^{(n,2)}(\mu_1) = -C_1 e^{-\beta \nu (2\pi n +\theta_2)} \sin \theta_2 - C_2 e^{-\beta (\nu+\mu_1/\beta)(2\pi n + \theta_2)} + e^{-\beta(2 \pi n+\theta_2)}.
\label{eq:asymphorns}
\end{array}	\right.
\end{equation}
On expanding $\sin \theta_{1}$ and $\sin \theta_2$ we get the expressions for two LP-branches forming the $n$-th `horn'
\begin{equation}
\hspace{-1cm}
\left\{
\begin{array}{ll}
	\mu_2^{(n,1)}(\mu_1) &= - e^{-\beta \nu (\theta^n_0 + \phi + 2\pi \delta_{i1})} \frac{e^{-2 \pi \beta \nu n}}{\sqrt{1+\beta^2\nu^2}} \left[ C_1 \left( 1-\frac{\beta^2\nu^2}{1+\beta^2\nu^2} \frac{C_2^2}{C_1^2}e^{-4\pi\mu_1 n} \right)^{1/2} \right. \\ & \hspace{1cm} \left.+ \frac{C_2}{\sqrt{1+\beta^2\nu^2}}e^{-\mu_1 (2 \pi n + 2 \pi \delta _{i1} +  \theta_0^n + \phi)} + O(\mu_1) \right].\\
 	\mu_2^{(n,2)}(\mu_1) &= - e^{-\beta \nu (\pi - \theta^n_0 + \phi)} \frac{e^{-2 \pi \beta \nu n}}{\sqrt{1+\beta^2\nu^2}} \left[- C_1 \left( 1-\frac{\beta^2\nu^2}{1+\beta^2\nu^2} \frac{C_2^2}{C_1^2}e^{-4\pi\mu_1 n} \right)^{1/2} \right. \\ & \hspace{1cm} \left.+ \frac{C_2}{\sqrt{1+\beta^2\nu^2}}e^{-\mu_1 (2 \pi n + \pi-\theta_0^n + \phi)} + O(\mu_1) \right].   
\end{array}   \right. 
\end{equation}
Upon setting $\mu_2^{(n)}$ to zero, we get a sequence $\left. \mu_1^{(n)} \right|_{\mu_2=0}$ of intersections of one of these branches with the axis $\mu_2 = 0$. Thus asymptotically
\begin{equation}\label{mu1zero}
\left. \mu_1^{(n)} \right|_{\mu_2=0} = \frac{1}{4\pi n}\left[\ln\left(\frac{C^2_2}{C^2_1}\right) + O\left(\frac{1}{n}\right)\right] .
\end{equation}
For genericity of the LP condition, we further require that the second derivative $\Phi_{\theta\theta}\neq 0$. Thus, the condition 
$$
\Phi_{\theta\theta}=0, 
$$
gives information about genericity and the location of the cusp point. We solve the following three conditions together
\begin{equation}
	\left\{\begin{array}{l}
	\Phi(\theta,\mu_1,\mu_2) = 0,\\  
	\Phi_\theta(\theta,\mu_1,\mu_2) = 0,\\
	\Phi_{\theta\theta}(\theta,\mu_1,\mu_2) = 0.
	\end{array}\right.
    \label{eq:cuspsys}
\end{equation}
Taking derivative with respect to $\theta$ in \eref{eq:trigterms} gives the third equation of \eref{eq:cuspsys},
\begin{equation}
 \cos(\theta-\phi)  +\frac{1}{\sqrt{1+\beta^2\nu^2}}\frac{C_2}{C_1}e^{-2\pi\mu_1 n}\left[\beta \nu \mu_1 + O(\mu_1^2)\right] + O(e^{-\alpha n}) = 0.
 \label{eq:cuspcond}
\end{equation}
Using \eref{eq:trigterms} and \eref{eq:cuspcond} we get
\begin{equation}
	\frac{1}{(1+\beta^2\nu^2)}\frac{C_2^2}{C_1^2}e^{-4\pi\mu_1 n}\left[\beta^2\nu^2+O(\mu_1)\right] + O(e^{-\alpha n}) = 1.
\end{equation}
which gives the value of $\mu_1$ at the cusp point,
\begin{equation}
\mu_1^n = \frac{1}{4\pi n} \left[\ln \left(   \frac{\beta^2\nu^2}{(1+\beta^2\nu^2)}\frac{C_2^2}{C_1^2}\right) +  O\left(\frac{1}{n}\right) \right].
	\label{eq:cuspfinal}
\end{equation}
The corresponding value of $\mu_2$ is obtained from \eref{eq:lpcond2}. We get,
\begin{equation}
 \mu_2^n = -e^{-\beta\nu(2\pi n + \theta_0 + \phi)} \frac{{\rm sign}(C_2)C_1}{\beta \nu \sqrt{1+\beta^2 \nu^2}} a^{-(\theta_0+\phi)/4\pi n} + O\left(\frac{1}{\sqrt{n}}\right),
\label{eq:cuspfinal2}
\end{equation} 
where $\theta_0$ is the value of $\theta^n_0$ at the cusp point, that is 
\begin{equation}
	\theta_0 = \left \{ \begin{array}{ll}
			\pi/2, & \textnormal{if} \ C_2<0,\\
			3\pi/2, & \textnormal{if} \ C_2>0,
			\end{array} \right.
\end{equation}
and
\begin{equation}
a = \frac{\beta^2\nu^2}{1+\beta^2\nu^2}\frac{C_2^2}{C_1^2}.
\end{equation}
Clearly, this cusp point is precisely where the two branches of a horn from \eref{eq:asymphorns} meet, i.e. when
$$\sin^2\theta^n_0= 1.$$
\subsubsection{PD curves}
The expressions derived to describe the LP-`horns' also describe PD bifurcation curves away from the cusp points. Indeed, let us now define
\begin{equation}
	\Psi(\theta,\mu_1,\mu_2) = \mu_2 + C_1 e^{-\beta \nu (2\pi n +\theta)} \sin \theta + C_2 e^{-\beta (\nu+\mu_1/\beta)(2\pi n + \theta)}.
\end{equation}
Then, the conditions for the PD-curve are
\begin{equation}
		\left\{\begin{array}{l}
		\Psi(\theta,\mu_1,\mu_2) = e^{-\beta(2 \pi n +\theta)},\\  
		\Psi_\theta(\theta,\mu_1,\mu_2) = \beta e^{-\beta(2\pi n + \theta)}.
		\end{array}\right.
	\label{eq:pdcond}
\end{equation}

On simplifying the second equation of \eref{eq:pdcond} we get 
\begin{equation}
\hspace{-1.5cm}
		\sin(\theta-\phi)  =-\frac{1}{\sqrt{1+\beta^2\nu^2}}\frac{C_2}{C_1}e^{-2\pi\mu_1 n}\left[\beta \nu - (1-\beta \nu \theta)\mu_1 + O(\mu_1^2)\right] - O( e^{-\alpha n}),
\end{equation}
which has the same leading terms as \eref{eq:trigterms}. We then follow the derivation for LP `horns' and obtain \eref{eq:thetapdlp} and \eref{eq:asymphorns} to describe PD curves. These asymptotics are only valid away from the cusp points, in a neighborhood of which we have to take into account higher order terms.\\

\subsection{Summarizing lemma for 1D model map}
We summarize our findings in the following lemma.

\begin{lemma}
\label{lem:horns}
In a neighborhood of the origin of the $(\mu_1,\mu_2)$-plane, the scalar model map $\eref{eq:1dmap0}$ has an infinite number of fold curves for fixed points $LP_{n}^{(1)}, n \in \mathbb{N},$ accumulating to the half axis $\mu_2=0$ with $\mu_1\geq0$. \\

Each curve resembles a `horn' with the following asymptotic representation of its two branches$:$
\begin{equation}
\hspace{-1.5cm}
\left\{
\begin{array}{ll}
	\mu_2^{(n,1)}(\mu_1) &= - e^{-\beta_0\nu_0 (\theta^n_0 + \phi_0 + 2\pi \delta_{i1})} \frac{e^{-2 \pi \beta_0 \nu_0 n}}{\sqrt{1+\beta_0^2\nu_0^2}} \left[ C^0_1 \left( 1-\frac{\beta_0^2\nu_0^2}{1+\beta_0^2\nu_0^2} \frac{(C^0_2)^2}{(C^0_1)^2}e^{-4\pi\mu_1 n} \right)^{1/2} \right. \\ & \hspace{1cm} \left.+ \frac{C^0_2}{\sqrt{1+\beta_0^2\nu_0^2}}e^{-\mu_1 (2 \pi n + 2 \pi \delta _{i1} +  \theta_0^n + \phi_0)} + O(\mu_1) \right].\\
 	\mu_2^{(n,2)}(\mu_1) &= - e^{-\beta_0 \nu_0 (\pi - \theta^n_0 + \phi_0)} \frac{e^{-2 \pi \beta_0 \nu_0 n}}{\sqrt{1+\beta_0^2\nu_0^2}} \left[- C^0_1 \left( 1-\frac{\beta_0^2\nu_0^2}{1+\beta_0^2\nu_0^2} \frac{(C^0_2)^2}{(C^0_1)^2}e^{-4\pi\mu_1 n} \right)^{1/2} \right. \\ & \hspace{1cm} \left.+ \frac{C^0_2}{\sqrt{1+\beta_0^2\nu_0^2}}e^{-\mu_1 (2 \pi n + \pi-\theta_0^n + \phi_0)} + O(\mu_1) \right].   
\end{array}   \right. 
\label{eq:horn_asymp}
\end{equation}
where 
\begin{equation*}
\begin{array}{l}
\phi_0 := \arcsin\left(\frac{1}{\sqrt{1+\beta_0^2 \nu_0^2}}\right),~~~
\theta^n_0 := \arcsin \left( -\frac{\beta_0 \nu_0}{\sqrt{1+\beta_0^2\nu_0^2}}\frac{C^0_2}{C^0_1}e^{-2\pi\mu_1 n}\right),
\end{array}
\end{equation*}
and $\delta_{ij}$ is the Kronecker delta where $i=-\textnormal{sign}(C^0_2)$.\\

The branches of each $LP_{n}^{(1)}$ curve meet at a cusp point $CP_{n}^{(1)}$ with the following asymptotic representation$:$
\begingroup
\renewcommand*{\arraystretch}{1.8}
\begin{equation}
\hspace{-1cm}
CP_n^{(1)}=
\left(
\begin{array}{l}
\mu_1^n\\
\mu_2^n
\end{array}\right) = 
\left(
\begin{array}{l}
 \frac{1}{4\pi n} \left[\ln \left(   \frac{\beta_0^2\nu_0^2}{(1+\beta_0^2\nu_0^2)}\frac{(C^0_2)^2}{(C^0_1)^2}\right) +  O\left(\frac{1}{n}\right) \right]\\
 -e^{-\beta_0\nu_0(2\pi n + \theta_0 + \phi_0)} \frac{{\rm sign}(C^0_2)C^0_1}{\beta_0 \nu_0 \sqrt{1+\beta_0^2 \nu_0^2}} a^{-(\theta_0+\phi_0)/4\pi n} + O\left(\frac{1}{\sqrt{n}}\right)
\end{array}\right),
\label{eq:cusps_asymp}
\end{equation}
\endgroup
where
\begin{equation*}
	\theta_0 := \left \{ \begin{array}{ll}
			\pi/2, & \textnormal{if} \ C^0_2<0,\\
			3\pi/2, & \textnormal{if} \ C^0_2>0,
			\end{array} \right.
\end{equation*}
and 
\begin{equation*}
a := \frac{\beta_0^2\nu_0^2}{1+\beta_0^2\nu_0^2}\frac{(C^0_2)^2}{(C^0_1)^2}\ .
\end{equation*}
\par
\medskip
Moreover, there exists an infinite number of period-doubling curves $PD_{n}^{(1)}, n \in \mathbb{N},$
which have -- away from the cusp points $CP_n^{(1)}$ -- the same asymptotic representation as the fold bifurcation curves $LP_n^{(1)}$. Depending on $(C^0_1,C^0_2)$, the period-doubling curves could either be smooth or have self-intersections developing small loops around the corresponding cusp points.
\end{lemma}



\section{Analysing the 3D model map}
\label{sec:3D_map}

In this section we study the original 3D model map \eref{eq:3Dmap} that we rewrite here for convenience as
\begin{eqnarray}
\hspace{-1.0cm}G: \left(\begin{array}{l}
x_1 \\ 
x_2 \\ 
x_4
\end{array}\right) \mapsto 	
\left(\begin{array}{l}
1 + \alpha_1 x_1 x_4^\nu \cos \left(-\frac{1}{\beta}\ln x_4 + \phi_1 \right) + \alpha_2 x_2 x_4^{\nu+\mu_1/\beta}\\
1 + \alpha_3 x_1 x_4^\nu \sin \left(-\frac{1}{\beta}\ln x_4 + \phi_2  \right) + \alpha_4 x_2 x_4^{\nu+\mu_1/\beta} \\
\mu_2 + C_1 x_1 x_4^\nu \sin \left(-\frac{1}{\beta}\ln x_4  \right) + C_2 x_2 x_4^{\nu+\mu_1/\beta}
\end{array}\right).
\label{eq:3dmap}
\end{eqnarray}
\par
\medskip
The analysis of fixed points of (\ref{eq:3dmap}) leads to the same equation (\ref{eq:1Dmap1}) for the $x_4$ coordinate. Thus, all conclusions about the existence and asymptotics of $LP_n^{(1))}$ and $PD_n^{(1)}$ curves, as well as $CP_n^{(1)}$ points in Lemma \ref{lem:horns}, remain valid. Indeed, taking into account the $o(|x|^{\nu})$-term does not alter the leading terms in any expression.

\subsection{Results of numerical continuation} 

\begin{figure}[ht!]
	\centering
	\begin{subfigure}[b]{0.45\textwidth}
		\includegraphics[width=\textwidth]{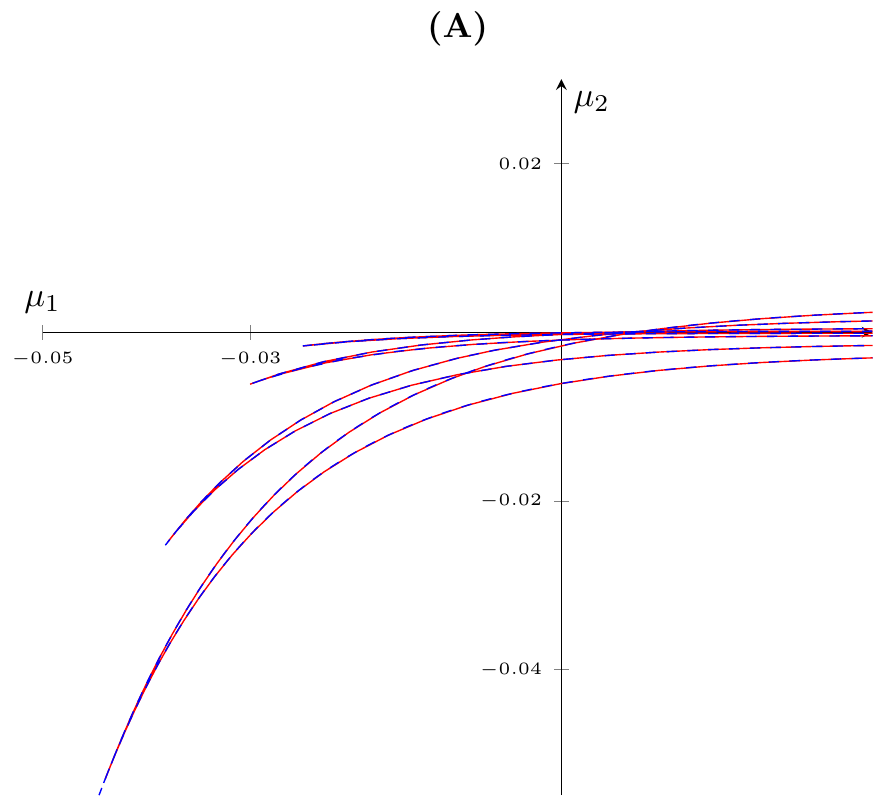}
	\end{subfigure}\hfill
	\begin{subfigure}[b]{0.48\textwidth}
		\includegraphics[width=\textwidth]{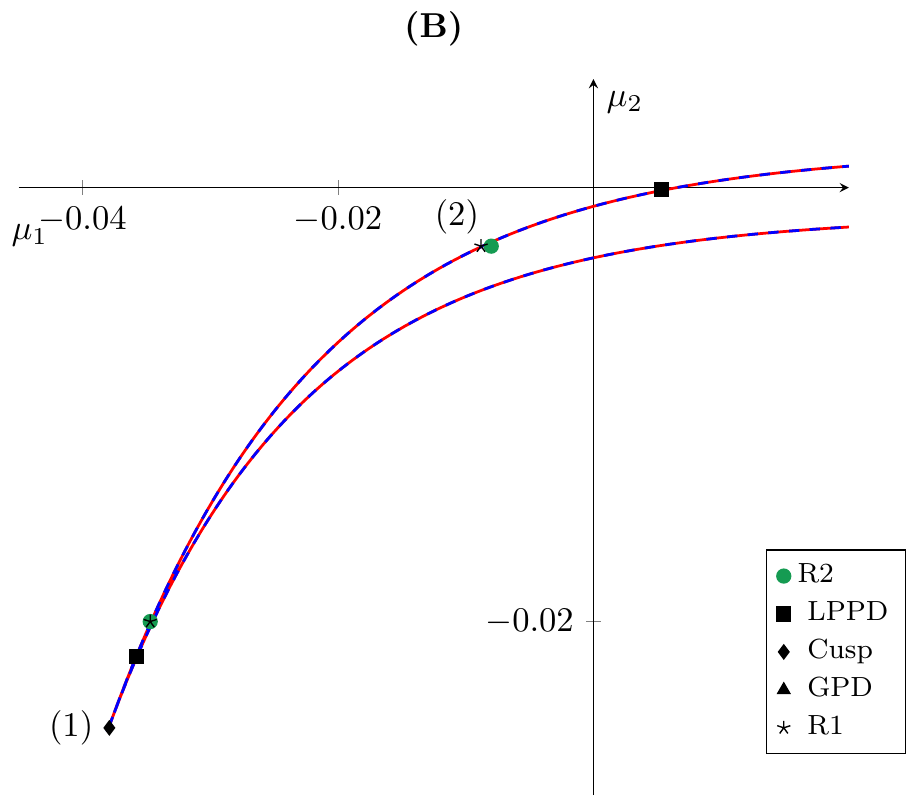}
	\end{subfigure}\\	
	\begin{subfigure}[b]{0.45\textwidth}
		\includegraphics[width=\textwidth]{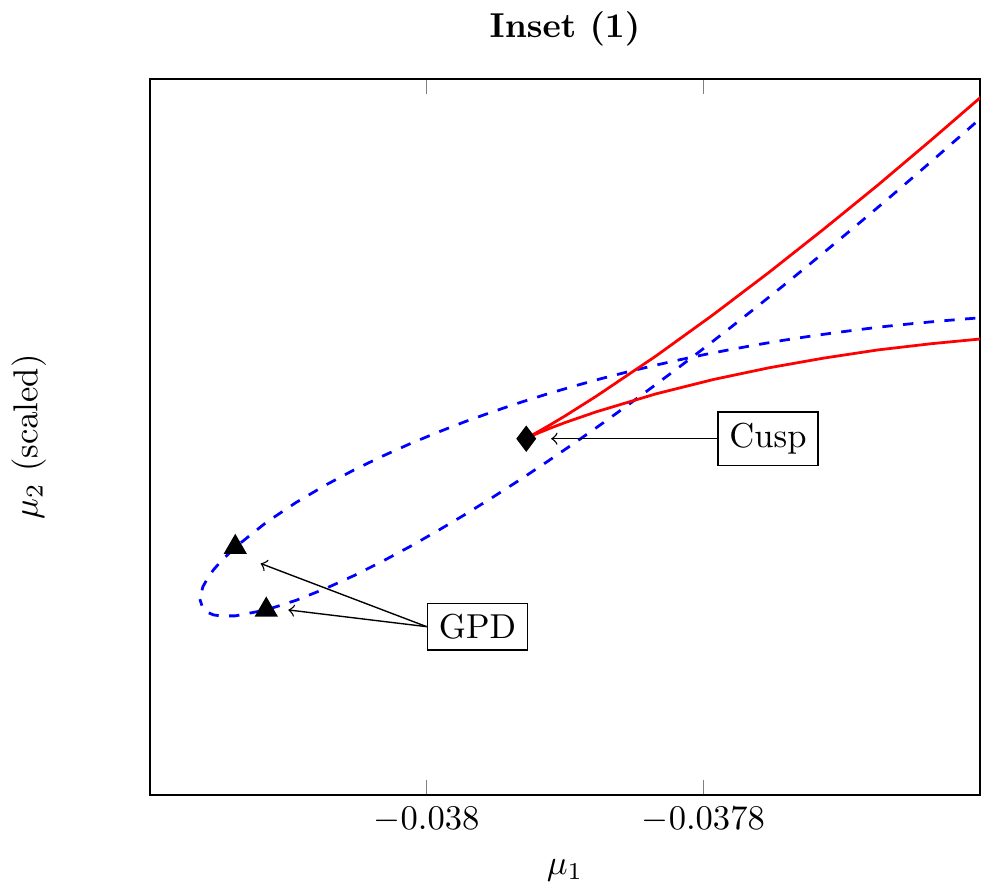}
	\end{subfigure}\hfill
	\begin{subfigure}[b]{0.45\textwidth}
		\includegraphics[width=\textwidth]{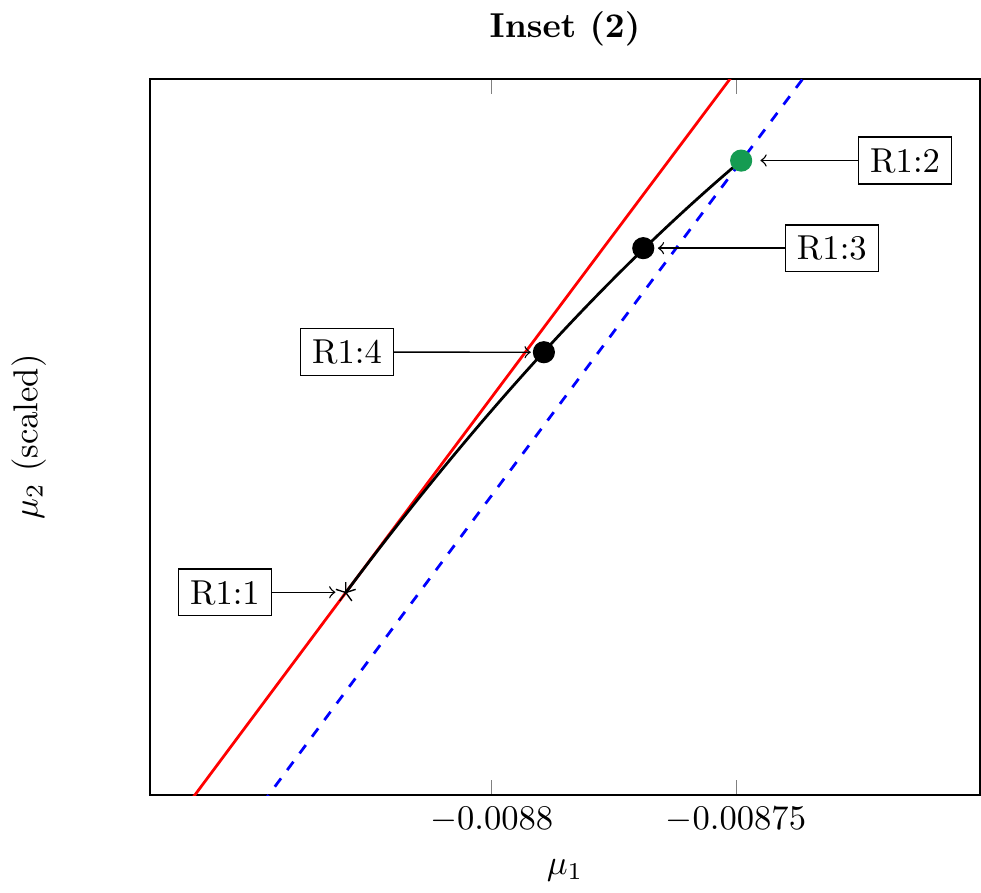}
	\end{subfigure}
	\caption{Primary LP (solid red) and PD (dashed blue)  curves obtained by numerical continuation for the map \eref{eq:3dmap}. The parameters are fixed as $\nu=0.5$, $\beta=0.5, C_1=0.8, C_2=1.2, \alpha_1 = 0.8, \alpha_2=1.3, \alpha_3 = 0.6, \alpha_4=1.1$ and $\phi_1=\phi_2 = \pi/6$. The curves have almost the same global structure as for the 1D map, as can be seen in \textbf{(A)}. In \textbf{(B)} one such curve is presented, together with several codim 2 points found along it. In  \textbf{Inset (1)} we see the previously described spring area made up by the PD and LP curves. Three codim 2 bifurcation points are observed, two corresponding to the generalised period doubling (GPD) bifurcation and one corresponding to the Cusp bifurcation. In \textbf{Inset (2)} we see the interaction between the 1:2 resonance (R2) point on the PD curve and the 1:1 resonance point (R1) on the LP curve, via the primary NS curve (solid black). On this curve we find two more codimension 2 bifurcation points: 1:3 resonance (R3) and 1:4 resonance (R4). }
	\label{fig:3Dcont}
\end{figure}

We look for fixed points of map (\ref{eq:3dmap}) and their various codim 1 curves. The results are similar to that of the scalar model map, except for higher dimensional codim 2 points that exist only in the 3D model map. In \autoref{fig:3Dcont}, we see the PD and LP curves obtained via numerical continuation for a fixed set of parameters: $\nu=\beta=0.5, C_1=0.8, C_2=1.2, \alpha_1 = 0.8, \alpha_2=1.3, \alpha_3 = 0.6, \alpha_4=1.1$ and $\phi_1=\phi_2 = \pi/6$. We immediately see similarities with the scalar case.
The global structure of these curves is the same as in the scalar case. They form sequences that accumulate onto the primary homoclinic curve asymptotically. The LP `horns' have cusp points and are accompanied by PD curves with/without GPD points (depending on saddle or spring area).
All this is expected as the scalar map is a correct asymptotic representation of the 3D model map. \\
There are however three main differences with respect to the scalar model map which can be attributed to the higher dimension of the 3D map:
\begin{enumerate}
    \item Spring and saddle areas may occur differently for the 1D and 3D model maps for the same parameter values. All these points numerically appear to asymptotically approach the 3DL transition point.
	\item Along the PD, LP and NS curves we observe many higher dimensional codimension 2 points. These points are R1:1 (Resonance 1:1), R1:2 (Resonance 1:2), LPPD (Fold-Flip), R1:3 (Resonance 1:3), R1:4 (Resonance 1:4).
	\item There exists two NS curves in a very small domain between the PD and LP curves. The end points of the NS segment are \textit{strong resonance} points.
\end{enumerate}
These points appear to numerically approach the origin $\mu=0$ (3DL transition). The endpoints of the NS curve are points R1 and R2, as can be seen in \autoref{fig:3Dcont} \textbf{(B)}. For a detailed discussion on the various codimension 2 points and their local bifurcation diagrams, see \cite{Kuznetsov2013}.

We did not see a significant difference in behavior of the PD/LP curves upon changing the coefficients $\alpha_i$ and $\phi_j$. This can be attributed to the effect of the corresponding terms in \eref{eq:3dmap} to the dynamics of $x_4$. These terms are  $o(\|x\|^{\nu})$ in the fixed point equation for $x_4$.\\
In \Tref{tab:codim2} we present sequences of some of the codimension 2 points found on successive PD/LP curves of \autoref{fig:3Dcont}. These sequences are obtained via detection along PD/LP curves from continuation. GPD and CP points are not reported as they are generally hard to detect along continuations, due to large test function values and absolute gradients. They are approximated in practice by noting where the sign of the scalar GPD test function changes. Note that codimension 2 points such as R1, R2 and LPPD were observed more than once on a single PD/LP curve. In \Tref{tab:codim2} we show only one point per curve for each of the different bifurcation points.

		\pgfplotstableset{
			begin table=\begin{longtable},
				end table={\label{tab:r2_1}\end{longtable}},
		}
		
		\begin{filecontents*}{datfiles/combined.dat}
   5.9031233e-03  -2.4502629e-04  -9.9024554e-03  -5.1621381e-03  -4.0228933e-02  -4.2268778e-02
   5.3053100e-03  -1.3752403e-04  -8.8301375e-03  -2.7386466e-03  -3.4664957e-02  -2.0046133e-02
   4.3915050e-03  -4.1359343e-05  -7.2678455e-03  -7.7386157e-04  -2.7399336e-02  -5.0254354e-03
   3.7374536e-03  -1.2119564e-05  -6.1795093e-03  -2.1924804e-04  -2.2749939e-02  -1.3355000e-03
   3.2506109e-03  -3.5142961e-06  -5.3758309e-03  -6.2199222e-05  -1.9480994e-02  -3.6375083e-04
   3.0515348e-03  -1.8883625e-06  -5.0477468e-03  -3.3139487e-05  -1.8180945e-02  -1.9083541e-04
   2.8753494e-03  -1.0137083e-06  -4.7574470e-03  -1.7659310e-05  -1.6045696e-02  -5.2879598e-05
   2.7183539e-03  -5.4376345e-07  -4.4987449e-03  -9.4114585e-06  -1.5157535e-02  -2.7905008e-05
   2.3356685e-03  -8.3650514e-08  -3.8678213e-03  -1.4254656e-06  -1.3648821e-02  -7.7985940e-06
\end{filecontents*}
\pgfplotstabletypeset[
		columns/0/.style={
			column name=$\mu_1$,
			dec sep align,
			/pgf/number format/fixed zerofill,
			/pgf/number format/precision=4,
		},
		columns/1/.style={
			column name=$\mu_2$,
			dec sep align,
			/pgf/number format/fixed zerofill,
			/pgf/number format/precision=4,
            column type/.add={}{|}
		},
		columns/2/.style={
			column name=$\mu_1$,
			dec sep align,
			/pgf/number format/fixed zerofill,
			/pgf/number format/precision=4
		},
		columns/3/.style={
			column name=$\mu_2$,
			dec sep align,
			/pgf/number format/fixed zerofill,
			/pgf/number format/precision=4,
            column type/.add={}{|}
		},
		columns/4/.style={
			column name=$\mu_1$,
			dec sep align,
			/pgf/number format/fixed zerofill,
			/pgf/number format/precision=4
		},
		columns/5/.style={
			column name=$\mu_2$,
			dec sep align,
			/pgf/number format/fixed zerofill,
			/pgf/number format/precision=4,
		},        
		every head row/.style={
			before row=\multicolumn{4}{c}{LPPD (1)}&
            \multicolumn{4}{c}{R1 (1)}&
            \multicolumn{4}{c}{R2 (2)}
            \\\toprule,    
			after row=\midrule
		},
		every last row/.style={
			after row=\bottomrule\\ \caption{Cascades of codimension two points numerically obtained during continuation of limit point/period doubling solutions of the 3D map \eref{eq:3dmap}. Note that the fixed parameter values are the same as in \autoref{fig:3Dcont}.}\label{tab:codim2}
		}]{datfiles/combined.dat}

\begin{figure}[ht!]
\includegraphics[width=\textwidth]{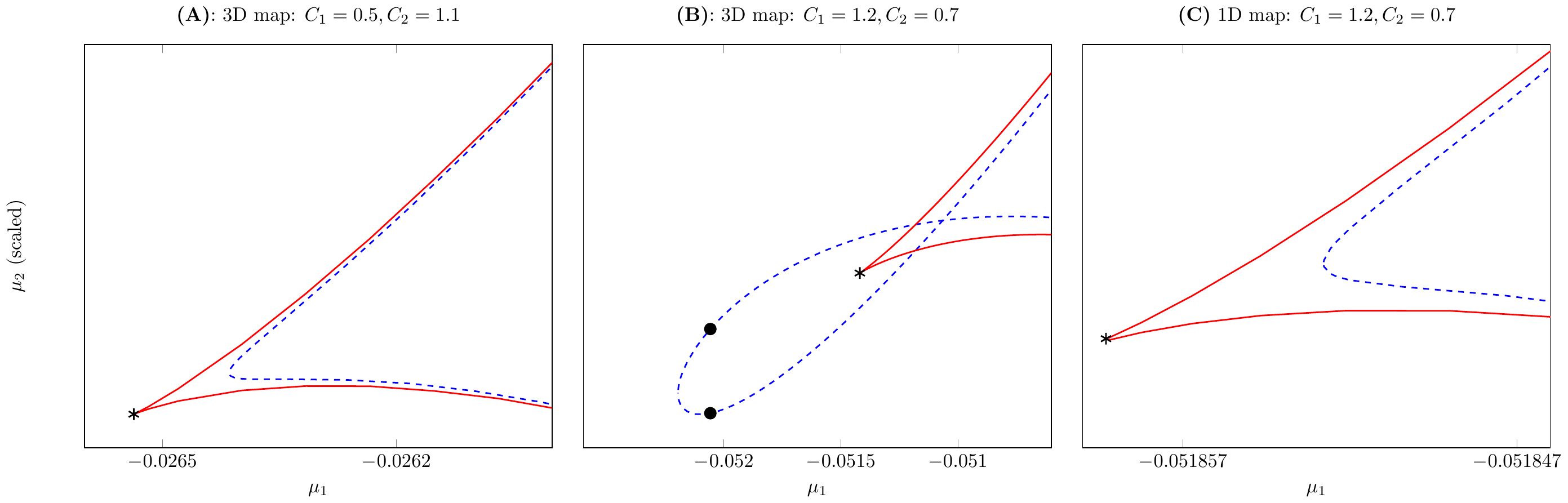}
	\caption{Plots of spring and saddle areas in the scalar map \eref{eq:1dmap} and 3D map \eref{eq:3Dmap}. In all plots $\mu_2$ is scaled for convenience. In \textbf{(A)} we see that there exists a saddle area in the 3D case, where GPD points are absent. \textbf{(B)} and \textbf{(C)} are plotted for the same value of $C_1$ and $C_2$, but with respect to the 3D map \eref{eq:3Dmap} and 1D map \eref{eq:1dmap} respectively. We see that the existence of the spring area in the 3D map does not imply the existence of the same in the 1D map. Other parameters fixed as in Figure \ref{fig:3Dcont}.}
	\label{fig:crossroad}
\end{figure}

For the scalar map we observed that there exist transitions between spring and saddle areas. These transitions can be explained by observing the appearance and disappearance of GPD points, as they exist generically on the PD loop in a spring area, and do not exist in the case of a saddle area. In the 3D case too, we numerically observe such transitions. However, when there is a spring (saddle) area in the 3D case, it does not imply that the same structure would exist in the 1D map for the same choice of parameters $C_1$ and $C_2$. This is shown in \autoref{fig:crossroad}. 

\subsection{Secondary homoclinic orbits}

In this section we analyze a particular type of homoclinic orbits, i.e. \textit{secondary} homoclinic orbits which -- after leaving the saddle along the unstable manifold -- make two global turns before returning to the saddle. \\

We look at the existence of these homoclinic orbits close to the primary homoclinic orbit in (\ref{eq:system}), upon perturbing parameters $\mu_1$ and $\mu_2$. The existence of the orbits is a codim 1 situation and corresponds to a curve in the $(\mu_1,\mu_2)$-plane. As before, we look for these curves in the wild case, where $\nu<1$. In the tame case $\nu>1$, they do not exist.\\

\begin{figure}[hbtp]
	\centering
	\includegraphics[width=0.45\textwidth]{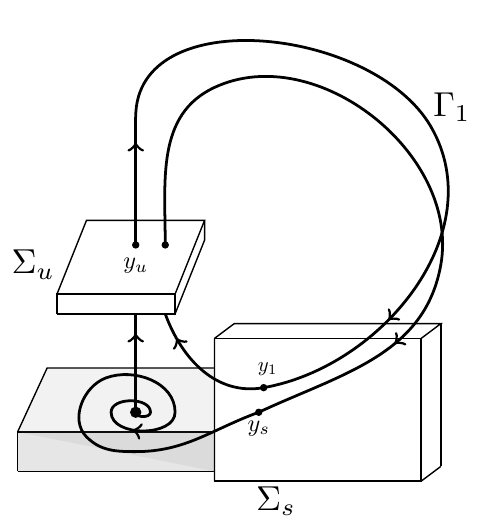}
	\caption{Poincar\'e map for the secondary homoclinic solution $\Gamma_1$. Upon leaving $y_u$ along the unstable manifold, the corresponding orbit makes two global turns and returns to the origin.}
	\label{fig:sechom_pmap}
\end{figure}

Consider \autoref{fig:sechom_pmap}. The secondary homoclinic orbit $\Gamma_1$ in the locally linearized ODE (\ref{eq:system}) leaves the point $y_u = (0,0,0,1) \in \Sigma_u$, along the unstable manifold and meets $\Sigma_s$ at $y_1=(1,0,1,\mu_2)$. From this point, the orbit departs again and this time returns along the stable manifold to approach the origin. The orbit crosses $\Sigma_s$ at $y_s=(1,0,1,0)$. Using the 3D model map $G$ defined by \eref{eq:3dmap}, the condition is
\begin{eqnarray}
G \left(\begin{array}{l}
1\\
1\\
\mu_2
\end{array}\right) = 
\left(\begin{array}{l}
1 \\ 
1\\
0
\end{array}\right),
\end{eqnarray}
which implies 
\begin{eqnarray}
\mu_2 + C_1 \mu_2^\nu \sin \left( -\frac{1}{\beta} \ln \mu_2 \right) + C_2 \mu_2^{\nu+\mu_1/\beta} = 0.
\label{eq:sechom}
\end{eqnarray}


Let us define 
\begin{eqnarray}
H(\mu) := \mu_2 + C_1 \mu_2^\nu \sin \left( -\frac{1}{\beta} \ln \mu_2 \right) + C_2 \mu_2^{\nu+\mu_1/\beta}.
\end{eqnarray}
Note that here $\mu_2$ must be positive. The shape of $H(\mu)=0$ is similar to the curve $F(x,\mu) = 0$ (from \eref{eq:1dmap}). For positive $\mu_1$, it is possible to obtain infinitely many solutions of \eref{eq:sechom} for $\mu_2$ sufficiently small. That is not the case when $\mu_1<0$, as there are only finitely many or no non-trivial solutions for $\mu_2$ sufficiently small.\\

\begin{figure}[hbtp]
	\centering
	\begin{subfigure}[b]{0.6\textwidth}
		\includegraphics[width=\textwidth]{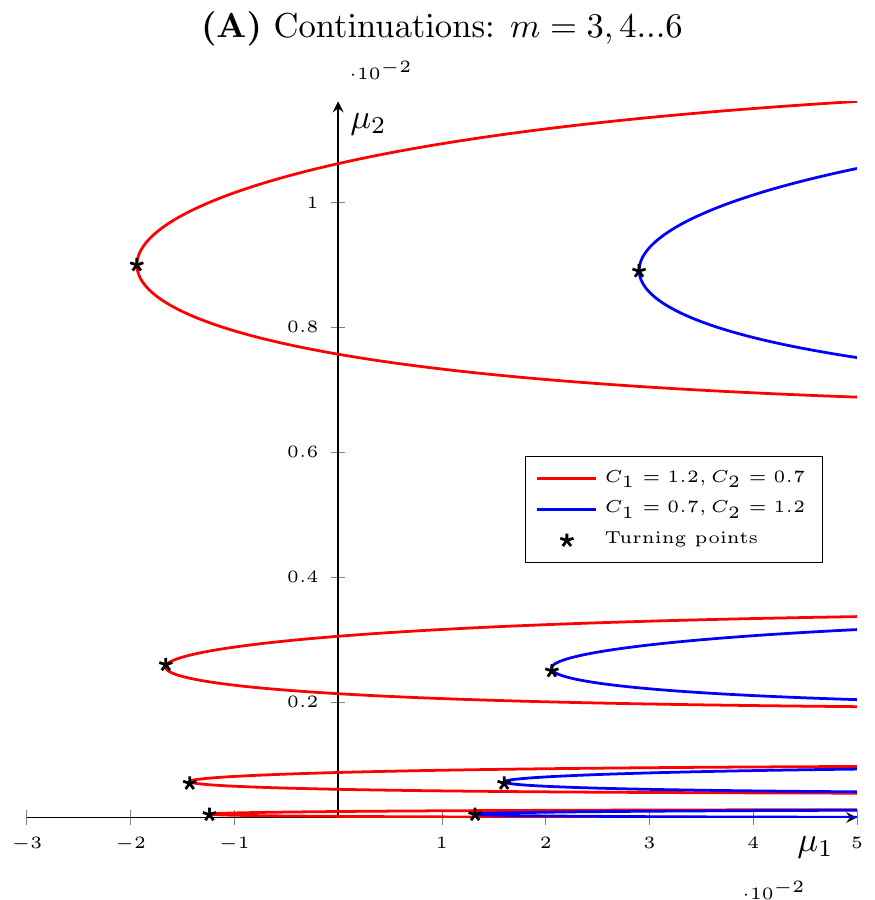}
	\end{subfigure}\\
	\begin{subfigure}[b]{0.5\textwidth}
		\includegraphics[width=\textwidth]{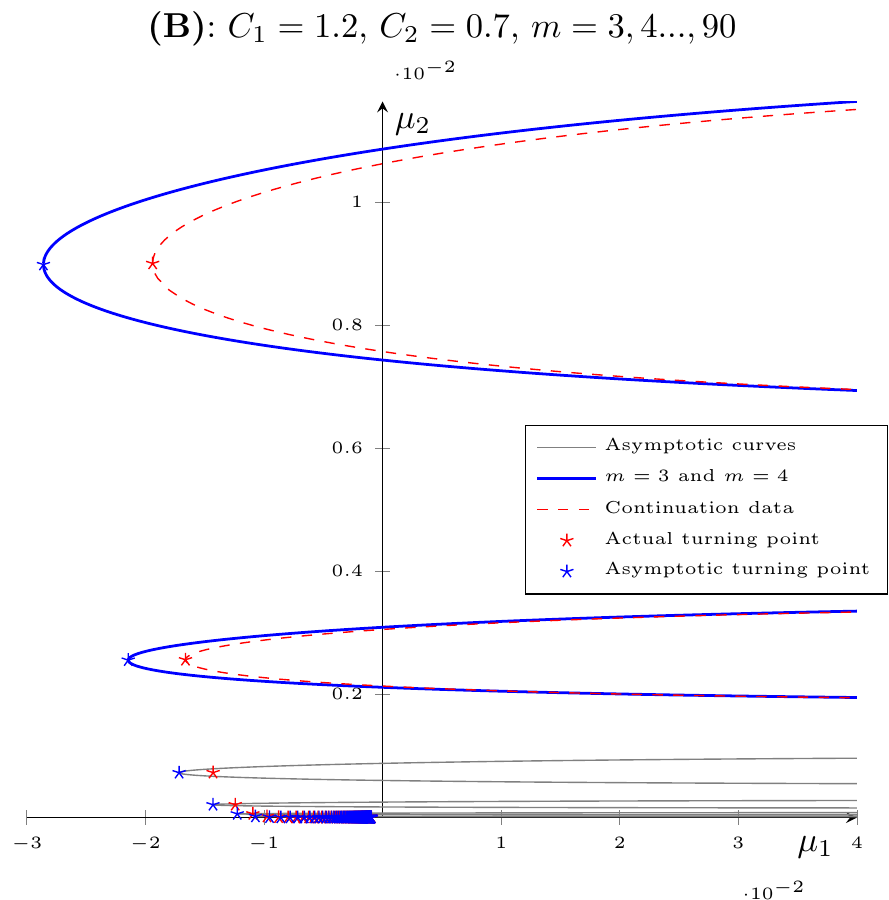}
	\end{subfigure}\hfill
	\begin{subfigure}[b]{0.45\textwidth}
		\includegraphics[width=\textwidth]{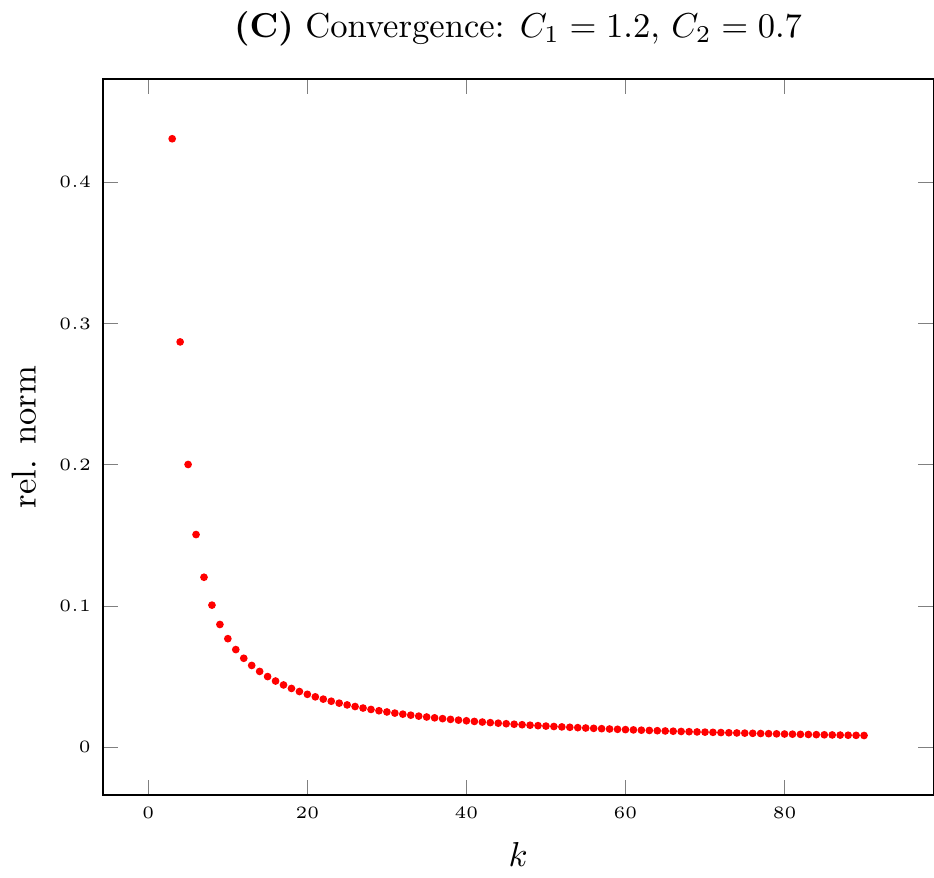}
	\end{subfigure}
	\caption{Solutions of \eref{eq:sechom} in $(\mu_1,\mu_2)$-space.  In \textbf{(A)}, `parabolas' are obtained via continuation in Matcont, for two sets of parameter values. The exact turning points (in black) are obtained by Netwon iterations. In \textbf{(B)}, exact curves are plotted along with asymptotic curves.  In \textbf{(C)}, we plot relative norm differences between asymptotic and exact turning points.}
	\label{fig:sechom}
\end{figure}

In \autoref{fig:sechom} the non-trivial solutions are continued with respect to the parameters $\mu_1$ and $\mu_2$ for two different sets of values of $C_1$ and $C_2$. We observe three things:
\begin{enumerate}
	\item There are secondary homoclinic curves which form \textit{horizontal parabolas} and these `parabolas' approach the primary homoclinic curve $\mu_2=0$ asymptotically.
	\item These `parabolas' possess \textit{turning} points where the two upper and lower secondary homoclinic branches merge. The sequence of turning points obtained from successive `parabolas' appears to approach the origin asymptotically.
	\item For different values of $C_1$ and $C_2$, the sequence of turning points is located strictly either in the first or second quadrant.
\end{enumerate}

\subsection{Asymptotics of secondary homoclinics}

The observations above can be explained to some extent by asymptotic expressions for the parabolas and the corresponding turning points.\\

\subsubsection{`Parabolas'}

Noticing $\mu_2>0$, let
\begin{eqnarray}
-\frac{1}{\beta} \ln \mu_2 = 2\pi m + \theta,
\label{eq:mu2parabola}
\end{eqnarray}
for large $m \in \mathbb{N}$  and $\theta \in (0, 2\pi)$.  On dividing both sides by $\mu_2^\nu \neq 0$ and using the above parametrization for $\mu_2$, \eref{eq:sechom} becomes
\begin{eqnarray}
e^{-\beta(1-\nu)(\pi m + \theta)} + C_1 \sin \theta + C_2 e^{-\mu_1 (\pi m +\theta)} = 0.
\label{eq:sechom2}
\end{eqnarray}
On simplifying, we get
\begin{equation}
\sin\theta = -\frac{C_2}{C_1}e^{-2\pi\mu_1 m}(1-\mu_1\theta+O(\mu_1^2)) + O(e^{-\alpha m}),
\label{eq:parcond1}
\end{equation}
where $\alpha = 2\pi\beta(1-\nu)$. For large $m$ and negative $\mu_1$, a solution $\theta$ exists only for small $|\mu_1|$. Thus we get two solutions $\theta$ from \eref{eq:parcond1},
\begin{equation}
\begin{array}{l}
\theta_1 = \theta^m_0+ \delta_{i1}(2\pi) + O(1/m),\\
\theta_2 = \pi - \theta^m_0 + O(1/m), 
\end{array}
\end{equation}
where 
$$
\theta^m_0 := \arcsin\left( -\frac{C_2}{C_1}e^{-2\pi\mu_1 m} \right),
$$
the index $i =-{\rm sign}(C_2)$ and $\delta_{ij}$ is the Kronecker delta. Thus the expressions for two `half-parabolas' are
\begin{equation}
\left\{
\begin{array}{ll}
\mu_2^{(m,1)} = e^{-\beta (2 \pi m + \theta^m_0 + \delta_{i1}(2\pi))}(1+O(1/m)),\\
\mu_2^{(m,2)} = e^{-\beta (2 \pi m + \pi - \theta^m_0 )}(1+O(1/m)).
\end{array}\right.
\end{equation}
Taking derivative with respect to $\theta$ in \eref{eq:parcond1} gives 
\begin{equation}
\cos\theta = \frac{C_2}{C_1}e^{-2\pi\mu_1 m}(\mu_1 + O(\mu_1^2)) + O(e^{-\alpha m}).
\label{eq:turncond}
\end{equation}
Solving \eref{eq:parcond1} and \eref{eq:turncond} together gives the condition for turning points. Using the two conditions gives,
\begin{equation}
\frac{C_2^2}{C_1^2}e^{-4\pi\mu_1 m}(1+O(\mu_1)) + O(e^{-\alpha m}) = 1.
\end{equation}
From this we get $\mu_1$,
\begin{equation}\label{mu1hom2}
\mu_1 = \frac{1}{4\pi m}\left[\ln \left( \frac{C_2^2}{C_1^2} \right) + O\left(\frac{1}{m}\right) \right],
\end{equation}
which is the same as the condition
\begin{equation}
\sin^2\theta = 1.
\end{equation}
Thus the  sequence of turning points is given by
\begingroup
\renewcommand*{\arraystretch}{1.3}
\begin{eqnarray}
\left(\begin{array}{l}
\mu_1^{(m)}\\
\mu_2^{(m)}
\end{array}\right) = 
\left(\begin{array}{l}
\frac{1}{(4\pi m)}\left ( \ln \left(\frac{C_2^2}{C_1^2}\right) + O\left(\frac{1}{m}\right)\right) \\
e^{-\beta(2\pi m + \theta_0)}\left(1+O\left(\frac{1}{m}\right)\right)
\end{array}\right),
\label{eq:sechomturning}
\end{eqnarray}
\endgroup
where 
\begin{equation}
	\theta_0 = \left \{ \begin{array}{ll}
			\pi/2, & \textnormal{if} \ C_2<0,\\
			3\pi/2, & \textnormal{if} \ C_2>0.
			\end{array} \right.
\end{equation}
We summarize the results in the following lemma.
\begin{lemma}
For the 3D model map $G$ defined by $\eref{eq:3Dmap}$, the condition 
\begin{eqnarray}
G \left(\begin{array}{l}
1\\
1\\
\mu_2
\end{array}\right) = 
\left(\begin{array}{l}
1 \\ 
1\\
0
\end{array}\right)
\end{eqnarray}
defines in a neighbourhood of the origin of the $(\mu_1,\mu_2)$-plane,  an infinite sequence of `parabolas' $Hom^{(2)}_m,~~ m \in \mathbb{N},$ that accumulate onto the half axis $\mu_2=0$ with $\mu_1 \geq 0$. Each parabola is formed by two branches with the following asymptotic representation:
\begin{equation}
\left\{
\begin{array}{ll}
\mu_2^{(m,1)} = e^{-\beta_0 (2 \pi m + \theta^m_0 + \delta_{i1}(2\pi))}\left(1+O\left(\frac{1}{m}\right)\right),\\
\mu_2^{(m,2)} = e^{-\beta_0 (2 \pi m + \pi - \theta^m_0 )}\left(1+O\left(\frac{1}{m}\right)\right),
\end{array}\right.
\label{eq:parab_asymp}
\end{equation}
where 
$$
\theta^m_0 := \arcsin\left( -\frac{C^0_2}{C^0_1}e^{-2\pi\mu_1 m} \right).
$$
\par
\medskip
These branches meet at a sequence of turning points $T^{(2)}_m$, which converges to the origin of the $(\mu_1,\mu_2)$-plane and is given by
\begin{equation}
T^{(2)}_m = 
\left(\begin{array}{l}
\mu_1^{(m)}\\
\mu_2^{(m)}
\end{array}\right) = 
\left(\begin{array}{l}
\frac{1}{(4\pi m)}\left ( \ln \left[\frac{(C^0_2)^2}{(C^0_1)^2}\right] + O\left(\frac{1}{m}\right)\right) \\
e^{-\beta_0(\pi m + \pi/2 + \delta_{i1}(2 \pi))}\left(1+O\left(\frac{1}{m}\right)\right)
\end{array}\right).
\label{eq:turning_asymp}
\end{equation}
\label{lem:parabola}
\end{lemma}

\section{Interpretation for the original ODE system}
\label{sec:ODE_interpretation}

Let us consider the original 4D system \eref{eq:system} in the linearizing coordinates near the equilibrium, the geometric construction in Figure \ref{fig:poincare_map} and the full 3D map $\Pi$ defined by \eref{eq:3DmapFull}. \\

Fixed points of this map $\Pi$ in $\Sigma_s$ correspond to periodic orbits, thus period doubling and fold bifurcations of these fixed points of this map correspond to the same bifurcations of periodic orbits in the original ODE system. \\


The second iterate of the map \eref{eq:3DmapFull}, such that $\mu_2>0$, defines an orbit in the original system \eref{eq:system} which makes an extra global excursion before returning to $\Sigma_u$. Starting at a point in the unstable 1D manifold of the equilibrium and letting the third component of the image to zero, implies that we consider an orbit of the ODE that departs along the unstable manifold and returns along the stable manifold back to the saddle. This orbit is therefore a secondary homoclinic orbit near the primary one. \\

Using Lemmas \ref{lem:horns} and \ref{lem:parabola} we are now able to formulate our main results in terms of the original 4D ODE near the 3DL-homoclinic transition. It follows from the fact that taking into account the $o(\|x\|^{\nu})$-term in \eref{eq:3DmapFull} does not alter the leading terms in all expressions, which further implies that given asymptotics are the same for the truncated map \eref{eq:3dmap} and full 3D return map \eref{eq:3DmapFull}.

\begin{theorem}
\label{Th:MAIN}
Consider a smooth 4D ODE system depending on two parameters
\begin{equation}
\label{4DODE}
\dot{x}=f(x,\alpha),~~x \in \R^4,~\alpha \in \R^2.
\end{equation}
Suppose that at $\alpha=0$ the system $(\ref{4DODE})$ satisfies the following assumptions:
\begin{enumerate}
	\item[] {\em \textbf{(B.1)}} The eigenvalues of the linearisation at the critical 3DL equilibrium $x=0$ are     
$$
\delta_0,\delta_0 \pm i \omega_0 \textnormal{ and } \epsilon_0,
$$  
where $\delta_0<0, \omega_0>0$, $\epsilon_0>0$ and $\sigma_0=\delta_0+\epsilon_0>0$.
	\item[] {\em \textbf{(B.2)}} There exists a primary homoclinic orbit $\Gamma_0$ to this 3DL equilibrium. 
	\item[] {\em \textbf{(B.3)}} The homoclinic orbit $\Gamma_0$ satisfies the following genericity condition: The normalized tangent vector to $\Gamma_0$ has nonzero projections to both the 1D eigenspace corresponding to the real eigenvalue $\delta_0$ and to the 2D eigenspace corresponding to the complex eigenvalues $\delta_0 \pm i \omega_0$, when approaching the equilibrium at $t\rightarrow \infty$.  	
\end{enumerate}

Then, in addition to the primary homoclinic curve $Hom^{(1)}$,  the bifurcation set of $(\ref{4DODE})$ in a neighborhood of $\alpha=0$ generically contains the following elements$:$
\begin{enumerate}
\item[{\em (i)}]
An infinite number of fold bifurcation curves $LP_{n}^{(1)}, n \in \mathbb{N},$ along which limit cycles with multiplier $+1$ exist making one global excursion and a number of small turns near the equilibrium. These curves accumulate to the saddle-focus part of the primary homoclinic curve. Each curve resembles a `horn' consisting of two branches that meet at a cusp point $CP_n^{(1)}$. The sequence of cusp points converges to $\alpha=0$.
\item[{\em (ii)}]
An infinite number of period-doubling bifurcation curves $PD_{n}^{(1)}, n \in \mathbb{N},$ along which limit cycles with multiplier $-1$ exist making one global excursion and a number of small turns near the equilibrium. Away from the cusp points $CP_n^{(1)}$, these period-doubling curves have the same asymptotic properties as the fold bifurcation curves $LP_n^{(1)}$. These period-doubling curves could either be smooth or have self-intersections developing small loops around the corresponding cusp points.
\item[{\em (iii)}]
An infinite number of secondary homoclinic curves $Hom_m^{(2)}, m \in \mathbb{N},$ along which the equilibrium has homoclinic orbits making two global excursions and a number of turns near the equilibrium after the first global excursion. These curves also accumulate to the saddle-focus part of the primary homoclinic curve. Each curve resembles a `parabola' and the sequence of turning points converges to $\alpha=0$.
\end{enumerate}
\end{theorem}

\begin{figure}[hbtp]
	\centering
	\includegraphics[width=\textwidth]{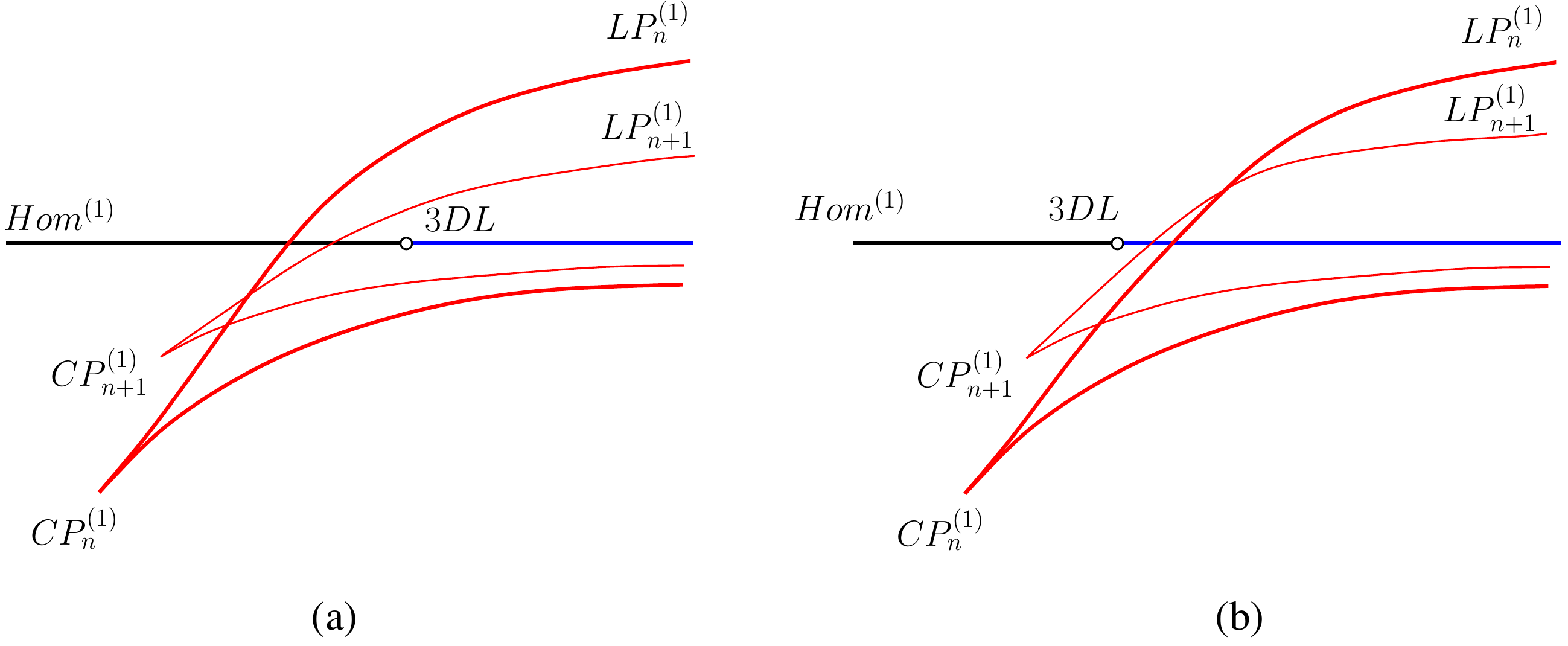}
	\caption{A sketch of two consecutive LP horns from Theorem \ref{Th:MAIN}. The saddle-focus part of $Hom^{(1)}$ branch is drawn in blue.}
	\label{fig:LPhorns}
\end{figure}

Part (i) of Theorem \ref{Th:MAIN} is illustrated in Figure \ref{fig:LPhorns}. Notice that $LP^{(1)}_n$ curves can intersect the primary homoclinic branch $Hom^{(1)}$ either at saddle points (Figure \ref{fig:LPhorns}(a)) or at saddle-focus points (Figure \ref{fig:LPhorns}(b)). In terms of the 1D (or 3D) model map these cases correspond to $C^0_2>0$ and $C^0_1>C^0_2$ or $C^0_1<C^0_2$, respectively. See equation \eref{mu1zero}.\\

\begin{figure}[hbtp]
	\centering
	\includegraphics[width=\textwidth]{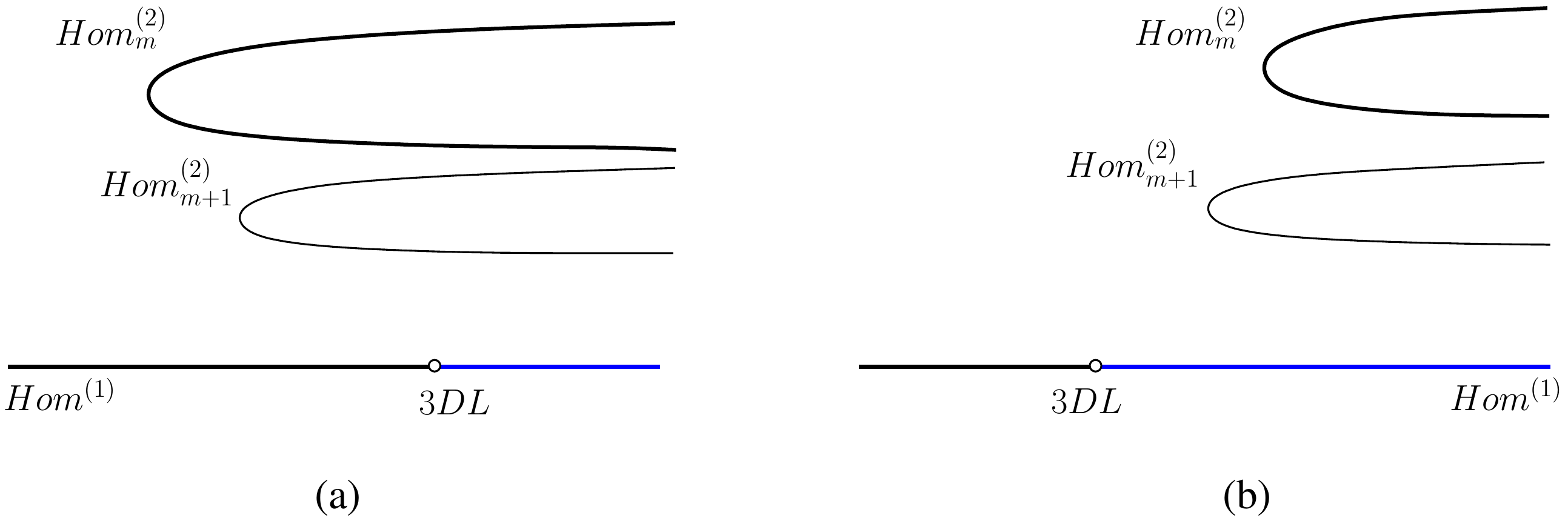}
	\caption{A sketch of two consecutive secondary homoclinic curves from Theorem \ref{Th:MAIN}. The saddle-focus part of $Hom^{(1)}$ branch is drawn in blue.}
	\label{fig:HOMsec}
\end{figure}

Part (iii) of Theorem \ref{Th:MAIN} is illustrated in Figure \ref{fig:HOMsec}. Notice that the turning points of the secondary homoclinic curves $Hom_m^{(2)}$ approach the 3DL transition point on $Hom^{(1)}$ either along its saddle part (Figure \ref{fig:HOMsec}(a)) or its saddle-focus part (Figure \ref{fig:HOMsec}(b)). In terms of the 2D (or 3D) model map these cases also correspond to $C^0_1>|C^0_2|$ or $C^0_1<|C^0_2|$, respectively. See equation \eref{mu1hom2}.\\

Our numerical analysis of the truncated model 3D map (\ref{eq:3Dmap}) also reveals 
NS curves in very small domains between the PD- and LP-curves. These curves correspond to {\em torus} bifurcation of cycles in the ODE system and do not exist for all combinations of $(C^0_1,C^0_2)$. The end points of the NS segment are strong resonance points. There are other codimension 2 points, i.e. GPD and LPPD. All these points should also be present in the generic ODE system and should form sequences that converge to the 3DL-transition point. 

\section{The Lorenz-Stenflo model}
\label{sec:example_model}

In this section we review the Lorenz-Stenflo model and discuss the presence of the 3DL-transition.\\

The Lorenz-Stenflo (LS) equations are a generalization of the well known Lorenz equations \cite{lorenz}, that describe low-frequency, short-wavelength acoustic-gravity perturbations in the atmosphere with additional dependence on the earth's rotation. The equations are as follows:
\begin{equation}
\label{Eq:LSmodel}
\left\{\begin{array}{rcl}
\dot{x} &=& \sigma(y-x) + su,\\
\dot{y} &=& rx - xz -y,\\
\dot{z} &=& xy - bz,\\
\dot{u} &=& -x - \sigma u,
\end{array}
\right.
\end{equation}
where $\sigma$ is the Prandtl number, $r$ is a generalized Rayleigh parameter, $b$ is a positive parameter and $s$ is a new parameter dependent on the Earth's rotation \cite{lorenzstenflo}. Setting $s=0$ reduces the first three equations in (\ref{Eq:LSperturbed}) back to the original Lorenz model.  \\

System (\ref{Eq:LSmodel}) possesses $\mathbb{Z}_2$-symmetry, and has one or three equilibria (the trivial equilibrium exists always). The system exhibits a wild 3DL transition, but the corresponding PD and LP curves are difficult to resolve due to highly contractive properties close to the transition, caused by large real parts of the eigenvalues at the trivial equilibrium.\\

To overcome this, we perturb the system to get:
\begin{equation}
\label{Eq:LSperturbed}
\left\{\begin{array}{rcl}
\dot{x} &=& \sigma(y-x) + su,\\
\dot{y} &=& rx - xz -y + \boldsymbol{\epsilon_1 z},\\
\dot{z} &=& xy - bz,\\
\dot{u} &=& -x - \sigma u + \boldsymbol{\epsilon_2 y} ,
\end{array}
\right.
\end{equation}
where the bold expressions are perturbation terms. This system is not $\mathbb{Z}_2$ symmetric anymore but still has a trivial equilibrium for all parameter values. We are not aware of any physical interpretation of the added terms.\\

\begin{figure}[hbtp]
\centering
\includegraphics[width=0.8\textwidth]{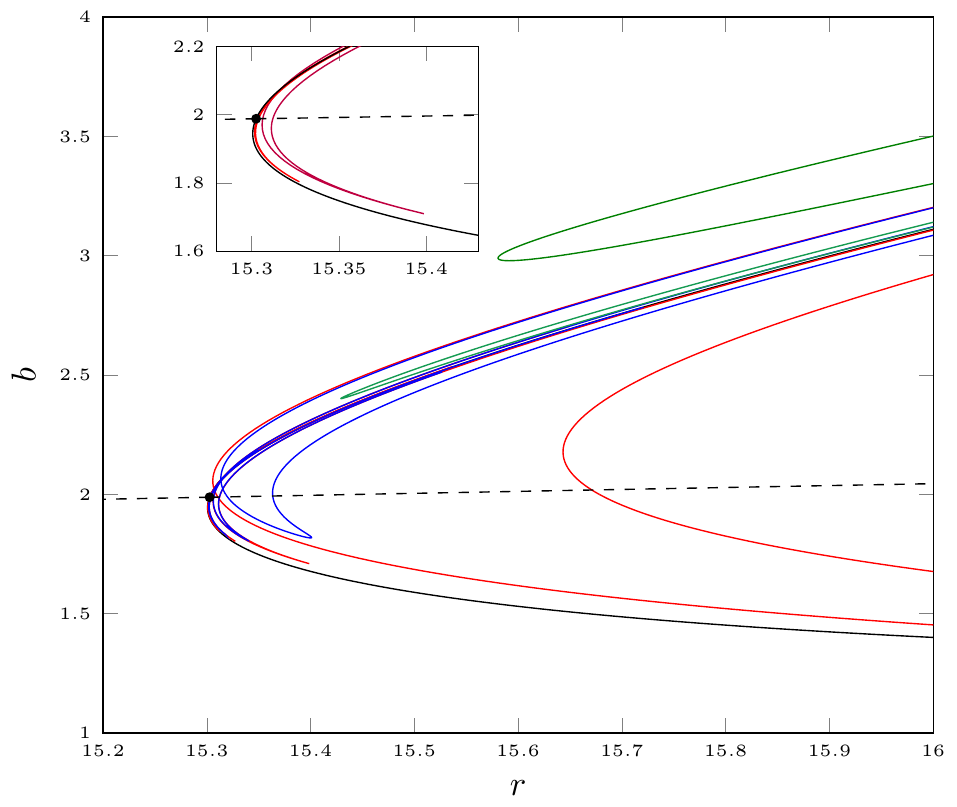}
\caption{Bifurcation curves near a wild 3DL transition in the $(b,r)$-plane: cyclic folds (red), period-doublings (blue), primary homoclinic (black), 3DL equilibrium transition (dashed black) and secondary homoclinics (green).}
\label{fig:3dl}
\end{figure}

In  \autoref{fig:3dl} we see a wild 3DL transition along the primary homoclinic curve (black) in the perturbed LS system (\ref{Eq:LSperturbed}) with
$$
\sigma=1,\ s=33,\ \epsilon_1=0.1, \ \epsilon_2=0.3.
$$
The 3DL-transition point is located at
$$
(r,b)\approx (15.302531 , 1.9884).
$$
The corresponding eigenvalues are $\delta_0 \pm i \omega_0, \delta_0$, and $\epsilon_0$ with  $\delta_0 \approx -1.9884 , \omega_0 \approx 6.2265, \epsilon_0 \approx 2.7769$.\\

We clearly see PD (blue) and LP (red) curves accumulating onto the primary homoclinic curve according to the theory. The PD curve within each `horn' forms a saddle area. The secondary homoclinic curves (green) form `parabolas' on one side of the primary homoclinic curve as expected. The curve of trivial equilibria with a 3D stable eigenspace is shown as a dashed line. The cusp points on each LP horn forms a sequence, and asymptotically approach the 3DL point at the intersection of the black curve with the dashed line. The inset shows only the LP horns. For this model the bifurcation curves have been computed using \textsc{matcont} \cite{Matcont} also based on \cite{ChaKuzSan1996,DeWitGovKuzFri2012}.\\

The 3DL-transition with the associated bifurcation structure exists also in the original $\mathbb{Z}_2$ symmetric system (\ref{Eq:LSmodel}). However, its bifurcation diagram will include additional bifurcation curves related to (symmetric) cycles and heteroclinic orbits. This will make the whole picture more complicated than in Figure \ref{fig:3dl}.

\section{Discussion}
\label{sec_discussion}
We have studied bifurcation diagrams of 4D two-parameter ODEs having at some critical parameter values a homoclinic orbit to an hyperbolic equilibrium with one simple unstable eigenvalue and three simple stable eigenvalues (one real and a complex-conjugate pair). We focused on the case, when a transition from a saddle homoclinic orbit to Shilnikov's wild saddle-focus homoclinic orbit takes place at the critical parameter values. Similar to the 3D Belyakov's saddle to wild saddle-focus homoclinic transition, we found infinite sequences of codim 1 bifurcations curves related to limit cycles (i.e., folds and period-doublings) and secondary homoclinic orbits accumulating on the primary (wild) saddle-focus homoclinic branch. However, there is a striking difference between these two cases. While in the standard Belyakov case all bifurcation curves approach the codim 2 point in the parameter plane tangentially to the saddle-focus homoclinic curve (having actually tangency of infinite order) and form `bunches', in the considered 3DL case neither of them approaches the codim 2 point. Instead, they form sequences of `horns' with cusps and other codim 2 bifurcation points, or `parabolas'. The sequences of codim 2 points and parabola tops indeed converge to the studied homoclinic 3DL point. In a sense, the bifurcation diagram for the considered 3DL-transition resembles more another codim 2 homoclinic bifurcation studied by Belyakov: A transition from tame to wild saddle-focus homoclinic orbit in 3D ODEs, when the saddle quantity vanishes \cite{Belyakov1984}. In that case, fold bifurcation curves for cycles also have cusp points accumulating to the transition point, while the secondary homoclinic curves look like `parabolas' with tops tending to the transition point. The exact source of this similarity remains a mystery but might be related to the simplicity of all eigenvalues in both cases.\\

The $C^1$-linearisation theorem used to derive the Poincar\'{e} return map (\ref{eq:3DmapFull}) permits to employ only the first-order partial derivatives of this map. This was sufficient to derive asymptotics for the fold and period-doubling bifurcations of the primary limit cycles, as well as those for the secondary homoclinic orbits. However, to verify nondegeneracy conditions for LP and PD bifurcations and to detect codim 2 points, one would need higher-order partial derivatives of the return map $\Pi$. Their existence can be guaranteed by using the $C^k$-linearisation  near the equilibrium with sufficiently big $k>1$. This exists generically, if one imposes a finite number of low-order non-resonance conditions on the \textit{real parts} of the eigenvalues of the equilibrium (see \cite{BroKop:1994}). However, it looks ptobable that one can avoid such extra conditions by careful analysis of the local flow near the equilibrium in the spirit of \cite[Sec. 2.8]{Shilnikovs1998}.\\

Clearly, 3DL-transitions are also possible in higher phase dimensions. The developed theory can be applied to such cases via a reduction to a four-dimensional homoclinic center-manifold \cite{Sandstede2000,ShaTur1999,Homburg1996} that exists near the bifurcation. However, also here some `gap' conditions should be imposed on the eigenvalues of the critical equilibrium to guarantee more than $C^1$ smoothness of the homoclinic center manifold. Whether one can avoid using the homoclinic center manifolds is also unclear.\\

Despite two above mentioned smoothness issues, we can conclude that, at least for a large open set of smooth ODEs, the 3DL homoclinic transition implies the complicated bifurcation structure described in the present paper.\\

An interesting challenge would be to prove analytically the existence of infinite sequences of generalized period-doubling points and strong resonances (see Figure \ref{fig:3Dcont}), at least for the truncated 3D model map (\ref{eq:3Dmap}).

\section*{References}
\bibliographystyle{unsrt}
\bibliography{references.bib}

\end{document}